\documentclass[11pt, twoside]{article}
\usepackage{amsmath,amsthm,amssymb,verbatim,enumerate,ifthen,times}
\usepackage[mathscr]{eucal}
\usepackage{enumerate}
\usepackage[utf8]{inputenc}
\def\titlerunning#1{\gdef\titrun{#1}}
\makeatletter
\def\author#1{\gdef\autrun{\def\and{\unskip, }#1}\gdef\@author{#1}}
\def\address#1{{\def\and{\\\hspace*{18pt}}\renewcommand{\thefootnote}{}%
\footnote {#1}}%
\markboth{\autrun}{\titrun}}
\makeatother
\def\email#1{e-mail: #1}
\def\subjclass#1{{\renewcommand{\thefootnote}{}%
\footnote{\emph{Mathematics Subject Classification (2010):} #1}}}
\def\keywords#1{\par\medskip
\noindent\textbf{Keywords.} #1}

\newcommand{\N}{\mathbb{N}}
\newcommand{\R}{\mathbb{R}}
\newcommand{\CC}{\mathbb{C}}
\newcommand{\C}{\mathscr{C}}
\newcommand{\CM}{\mathscr{CM}}
\newcommand{\CP}{\mathscr{CP}}
\newcommand{\A}{\mathscr{A}}
\newcommand{\B}{\mathscr{B}}

\newcommand{\D}{\mathscr{D}}

\def\d{\,{\rm d}}

\newtheorem{theorem}{Theorem}[section]
\newtheorem*{theorem*}{Theorem}
\def\Thm#1#2{\ifthenelse{\equal{#1}{*}}{\begin{theorem*}#2\end{theorem*}}
  {\begin{theorem}\label{T#1}#2\end{theorem}}}
\newtheorem{Atheorem}{Theorem}

\def\thm#1{Theorem~\ref{T#1}}

\newtheorem{proposition}[theorem]{Proposition}
\newtheorem*{proposition*}{Proposition}
\def\Prp#1#2{\ifthenelse{\equal{#1}{*}}{\begin{proposition*}#2\end{proposition*}}
             {\begin{proposition}\label{P#1}#2\end{proposition}}}

\newtheorem{corollary}[theorem]{Corollary}
\newtheorem*{corollary*}{Corollary}
\def\Cor#1#2{\ifthenelse{\equal{#1}{*}}{\begin{corollary*}#2\end{corollary*}}
             {\begin{corollary}\label{C#1}#2\end{corollary}}}
\def\cor#1{Corollary~\ref{C#1}}

\newtheorem{lemma}[theorem]{Lemma}
\newtheorem*{lemma*}{Lemma}
\def\Lem#1#2{\ifthenelse{\equal{#1}{*}}{\begin{lemma*}#2\end{lemma*}}
             {\begin{lemma}\label{L#1}#2\end{lemma}}}
\def\lem#1{Lemma~\ref{L#1}}

\newtheorem{example}[theorem]{Example}
\newtheorem*{example*}{Example}
\def\Exa#1#2{\ifthenelse{\equal{#1}{*}}{\begin{example*}\rm #2\end{example*}}
             {\begin{example}\label{Ex#1}\rm #2\end{example}}}

\newtheorem{problem}[theorem]{Problem}

\theoremstyle{definition}
\newtheorem{definition}[theorem]{Definition}

\newtheorem{remark}[theorem]{Remark}
\newtheorem*{remark*}{Remark}
\def\Rem#1#2{\ifthenelse{\equal{#1}{*}}{\begin{remark*}\rm #2\end{remark*}}
             {\begin{remark}\label{R#1}\rm #2\end{remark}}}

\newcommand{\eq}[1]{\eqref{E#1}}
\newcommand{\Eq}[2]{\ifthenelse{\equal{#1}{*}}
  {\begin{equation*}\begin{aligned}[]#2\end{aligned}\end{equation*}}
  {\begin{equation}\begin{aligned}[]\label{E#1}#2\end{aligned}\end{equation}}}
\newcommand{\DET}[1]{\begin{vmatrix}#1\end{vmatrix}}

\long\def\comment#1{}

\begin{document}


\baselineskip=17pt


\titlerunning{Equality and homogeneity of generalized integral means}

\title{Equality and homogeneity of generalized integral means}

\author{Zsolt P\'ales
\and 
Amr Zakaria}

\date{\today}

\maketitle

\address{Zs. P\'ales: Institute of Mathematics, University of Debrecen, H-4002 Debrecen, Pf.\ 400, Hungary;\\ 
\email{pales@science.unideb.hu}
\and
A. Zakaria: Department of Mathematics, Faculty of Education, Ain Shams University, Cairo 11341, Egypt;\\ 
\email{amr.zakaria@edu.asu.edu.eg}}
\subjclass{Primary 26D10, 26D15; Secondary 26B25, 39B72, 41A50\\This research has been supported by the Hungarian
Scientific Research Fund (OKTA) Grant K-111651.}


\begin{abstract}
Given two continuous functions $f,g:I\to\R$ such that $g$ is positive and $f/g$ is strictly monotone, a
measurable space $(T,\A)$, a measurable family of $d$-variable means $m: I^d\times T\to I$, and a probability
measure $\mu$ on the measurable sets $\A$, the $d$-variable mean $M_{f,g,m;\mu}:I^d\to I$ is defined by
$$
   M_{f,g,m;\mu}(\pmb{x})
      :=\left(\frac{f}{g}\right)^{-1}\left(
            \frac{\int_T f\big(m(x_1,\dots,x_d,t)\big)\d\mu(t)}
                 {\int_T g\big(m(x_1,\dots,x_d,t)\big)\d\mu(t)}\right)
   \qquad(\pmb{x}=(x_1,\dots,x_d)\in I^d).
$$
The aim of this paper is to solve the equality and homogeneity problems of these means, i.e., to find
conditions for the generating functions $(f,g)$ and $(h,k)$, for the family of means $m$, and for
the measure $\mu$ such that the equality
$$
   M_{f,g,m;\mu}(\pmb{x})=M_{h,k,m;\mu}(\pmb{x}) \qquad(\pmb{x}\in I^d)
$$
and the homogeneity property
$$
   M_{f,g,m;\mu}(\lambda\pmb{x})=\lambda M_{f,g,m;\mu}(\pmb{x}) 
   \qquad(\lambda>0,\,\pmb{x},\lambda\pmb{x}\in I^d),
$$
respectively, be satisfied.

\end{abstract}

\section{Introduction}
Throughout this paper, the symbols $\N$, $\R$, and $\R_+$ will stand for the sets of natural,
real, and positive real numbers, respectively, and $I$ will always denote a nonempty open real interval. 

In the sequel, a function $M:I^d\to I$ is called a
\emph{$d$-variable mean} on $I$ if the following so-called mean value property
\Eq{1}{
  \min(x_1,\dots,x_d)\leq M(\pmb{x})\leq \max(x_1,\dots,x_d)   \qquad(\pmb{x}=(x_1,\dots,x_d)\in I^d)
}
holds. Also, if both of the inequalities in \eq{1} are strict for all $x_1,\dots,x_d\in I$ with $x_i\neq x_j$
for some $i\neq j$, then we say that $M$ is a \emph{strict mean} on $I$. 

For various classes of means, the comparison, equality and homogeneity problems are of great importance. In 
this paper we aim to focus on the last two problems. In more details, the \emph{equality problem} in a class 
of $d$-variable means $\mathscr{M}_d(I)$ (defined on an interval $I$) of means is to find necessary and 
sufficient conditions in order that, for some means $M,N\in\mathscr{M}_d(I)$, the equality
\Eq{*}{
  M(\pmb{x})=N(\pmb{x})\qquad(\pmb{x}\in I^d)
}
be valid. For the formulation of the \emph{homogeneity problem}, we have to assume that $I\subseteq\R_+$. 
Then this problem is to determine all means $M\in\mathscr{M}_d(I)$ satisfying the functional equation
\Eq{}{
  M(\lambda \pmb{x})=\lambda M(\pmb{x})
  \qquad (\lambda>0,\,\pmb{x},\lambda \pmb{x}\in I^d).
}

In what follows, we describe many important classes of means and the solutions of the equality and 
homogeneity problems related to them. 

If $p$ is a real number, then the $d$-variable \textit{H\"older mean} $H_p:\R_+^d\to\R$ is defined as
\Eq{*}{
  H_{p}(\pmb{x})
  :=\begin{cases}
    \left(\dfrac{x_1^p+\cdots+x^p_d}{d}\right)^{\frac{1}{p}} &\mbox{if } p\neq0\\[3mm]
    \sqrt[d]{x_1\cdots x_d} &\mbox{if } p=0
   \end{cases}
  \qquad\big(\pmb{x}=(x_1,\dots,x_d)\in\R_+^d\big).
}
Obviously, $H_1$ and $H_0$ equal the arithmetic and geometric mean, respectively. It is easy to see that
H\"older means are strict and homogeneous means. 

In the sequel, the classes of continuous strictly monotone and continuous positive real-valued functions 
defined on $I$ will be denoted by $\CM(I)$ and $\CP(I)$, respectively. A classical generalization of H\"older 
means is the notion of \emph{$d$-variable quasi-arithmetic mean} (cf.\ \cite{HarLitPol34}), which is 
introduced as follows: For $f\in\CM(I)$ define
\Eq{A}{
  A_f(\pmb{x}):=f^{-1}\left(\frac{f(x_1)+\cdots+f(x_d)}{d}\right)
  \qquad\big(\pmb{x}=(x_1,\dots,x_d)\in I^d\big).
}
More generally, if $S_d$ denotes the $(d-1)$-dimensional simplex given by
\Eq{S}{
  S_d:=\{(t_1,\dots,t_d)\mid t_1,\dots,t_d\geq0,\,t_1+\dots+t_d=1\},
}
then we can also define
\Eq{At}{
  A_f(\pmb{x},\pmb{t}):=f^{-1}\big(t_1f(x_1)+\cdots+t_df(x_d)\big)
  \qquad\big(\pmb{x}=(x_1,\dots,x_d)\in I^d,\, \pmb{t}=(t_1,\dots,t_d)\in S_d\big),
}
which is called the \emph{weighted $d$-variable quasi-arithmetic mean on $I$.}
It is a classical result that, given two functions $f,g\in\CM(I)$, the equality of the generated 
quasi-arithmetic means $A_f$ and $A_g$ is equivalent to the existence of real constants $a,b$ such that 
$g=af+b$. Concerning the homogeneity problem related to quasi-arithmetic means, it is well-known that a
quasi-arithmetic mean is homogeneous if and only if it is a H\"older mean.

In this paper, we consider a much more general class of means. For their definition, we recall the notion of
Chebyshev system. We say that a pair $(f,g)$ of continuous functions defined on $I$ forms a
\textit{(two-dimensional) Chebyshev system on $I$} if, for any distinct elements $x,y$ of $I$, the
determinant
\Eq{*}{
  \D_{f,g}(x,y):=\DET{f(x) & f(y) \\ g(x) & g(y)}\qquad(x,y\in I)
}
is different from zero. If, for $x<y$, this determinant is positive, then $(f,g)$ is called a \emph{positive
system}, otherwise we call $(f,g)$ a \emph{negative system}. Due to the connectedness of the triangle
$\{(x,y)\mid x<y,\,x,y\in I\}$, it follows that every Chebyshev system is either positive or negative.
Obviously, if $(f,g)$ is a positive Chebyshev system, then $(g,f)$ is a negative one.

The most standard positive Chebyshev system on $\R$ is given by $f(x)=1$ and
$g(x)=x$. More generally, if $f,g:I\to\R$ are continuous functions with
$g\in\CP(I)$, $f/g\in\CM(I)$, then $(f,g)$ is a Chebyshev system. Indeed, we have
\Eq{D}{
  \D_{f,g}(x,y):=\DET{f(x)&f(y)\\g(x)&g(y)}
     =g(x)g(y)\left(\frac{f(x)}{g(x)}-\frac{f(y)}{g(y)}\right)
  \qquad(x,y\in I).
}
From, here it is obvious that $\D_{f,g}(x,y)$ vanishes if and only if $x=y$. Moreover, if $f/g$ is
decreasing (resp.\ increasing), then, for $x<y$, we have that $\D_{f,g}(x,y)>0$ (resp.\ $\D_{f,g}(x,y)<0$),
i.e., $(f,g)$ is a positive (resp.\ negative) Chebyshev system. By symmetry, analogous properties can be
established if $f$ is positive and $g/f$ strictly monotone.

For the sake of convenience and brevity, now we make the following hypotheses. We say that $m:I^d\times T\to
I$ is a \emph{measurable family of $d$-variable means on $I$} if
\begin{enumerate}[(H1)]
 \item $I$ is a nonvoid open real interval,
 \item $(T,\A)$ is a measurable space, where $\A$ is the $\sigma$-algebra of measurable sets of $T$,
 \item for all $t\in T$, $m(\cdot,t)$ is a $d$-variable mean on $I$,
 \item for all $\pmb{x}\in I^d$, the function $m(\pmb{x},\cdot)$ is measurable over $T$.
\end{enumerate}
If, instead of (H2) and H(4), we have that
\begin{enumerate}
 \item[(H2+)] $T$ is a topological space and $\A$ equals the $\sigma$-algebra $\B(T)$
of the Borel sets of $T$,
 \item[(H4+)] for all $\pmb{x}\in I^d$, the function $m(\pmb{x},\cdot)$ is continuous over $T$,
\end{enumerate}
then $m:I^d\times T\to I$ will be called a \emph{continuous family of $d$-variable means on $I$}.

The following lemma is the key to construct a mean in terms of a Chebyshev system, a measurable family of
means, and a probability measure (cf.\ \cite{PalZak17}).

\Lem{0}{Let $m:I^d\times T\to I$ be a measurable family of $d$-variable means, let $\mu$ be a probability
measure on $(T,\A)$ and let $(f,g)$ be a Chebyshev system on $I$. Then, for all $\pmb{x}\in I^d$, there exists
a unique element $y\in I$ such that
\Eq{Eq}{
  \int_T \D_{f,g}(m(\pmb{x},t),y)\d\mu(t)=0.
}
In addition, if $g$ is positive and $f/g$ is strictly monotone, then
\Eq{*}{
   y=\left(\frac{f}{g}\right)^{-1}\left(
            \frac{\int_T f\big(m(\pmb{x},t)\big)\d\mu(t)}
                 {\int_T g\big(m(\pmb{x},t)\big)\d\mu(t)}\right).
}}

The above lemma allows us to define a $d$-variable mean $M_{f,g,m;\mu}:I^d\to I$. Given $\pmb{x}\in I^d$, let
$M_{f,g,m;\mu}(\pmb{x})$ denote the unique solution $y$ of equation \eq{Eq}. In the particular case when $g$
is positive and $f/g$ is strictly monotone, we have that
\Eq{fg}{
   M_{f,g,m;\mu}(\pmb{x})
      :=\left(\frac{f}{g}\right)^{-1}\left(
            \frac{\int_T f\big(m(\pmb{x},t)\big)\d\mu(t)}
                 {\int_T g\big(m(\pmb{x},t)\big)\d\mu(t)}\right)
   \qquad(\pmb{x}\in I^d).
}
This mean will be called a $d$-variable \textit{generalized Bajraktarevi\'c mean} in the sequel.
In the case when $m$ is a two-variable family of weighted arithmetic means, this class of means was 
introduced and their comparison problem was also solved in the paper \cite{LosPal08}.   When $g=1$, then
\Eq{*}{
 M_{f,1,m;\mu}(\pmb{x})=(f)^{-1}\left(\int_T f\big(m(\pmb{x},t)\big)\d\mu(t)\right)
   \qquad(\pmb{x}\in I^d)
}
which will be termed a $d$-variable \textit{generalized quasi-arithmetic mean}.
If $T=\{1,\dots,d\}$, $\mu=\frac{\delta_1+\dots+\delta_d}{d}$ (where $\delta_{t}$ denotes the Dirac measure
concentrated at $t$) and $m(\pmb{x},t)=x_t$, then
\Eq{*}{
   M_{f,g,m;\mu}(\pmb{x})=M_{f,g}(\pmb{x})
    :=\left(\frac{f}{g}\right)^{-1}\bigg(\frac{f(x_1)+\dots+f(x_d)}{g(x_1)+\dots+g(x_d)}\bigg)
    \qquad\big(\pmb{x}=(x_1\dots,x_d)\in I^d\big),
}
which was introduced and studied by Bajraktarevi\'c \cite{Baj58}, \cite{Baj69}. In these two papers the 
equality problem of (standard) Bajraktarevi\'c means was solved under two times differentiability assumptions. 
The homogeneity of these means (assuming it for all $d\geq2$) was characterized in terms of Gini means (see 
their definition below) by Acz\'el and Dar\'oczy in \cite{AczDar63c}. When $g=1$, then 
$M_{f,1,m;\mu}(\pmb{x})=A_{f}(\pmb{x})$, which is the $d$-variable quasi-arithmetic mean introduced in \eq{A}.

To define the $d$-variable \textit{generalized Gini means}, let $p,q\in\CC$ such that either $p,q\in\R$ or 
$p=\bar{q}$ (this holds if and only if $p+q$ and $pq$ are real numbers or, equivalently, if $p$ and $q$ are 
the roots of a second degree polynomial with real coefficients). For $x\in\R_+$ define
\Eq{pq}{
  f(x)&:=x^p,&\quad g(x)&:=x^q& \quad&\mbox{if}\quad p,q\in\R,\,p\neq q,\\
  f(x)&:=x^p\log(x),&\quad g(x)&:=x^p& \quad&\mbox{if}\quad p=q\in\R,\\
  f(x)&:=x^a\sin(b\log(x)),&\quad g(x)&:=x^a\cos(b\log(x))&\quad&\mbox{if}\quad
      p=\bar{q}=a+bi,\,\,a,b\in\R,\,b\neq0.
}
After a simple computation, we have that
\Eq{*}{
  \D_{f,g}(x,y)&=x^{p+q}\big(\big(\tfrac{y}{x}\big)^q-\big(\tfrac{y}{x}\big)^p\big)
   \qquad&&\mbox{if}\quad p,q\in\R,\,p\neq q,\\
  \D_{f,g}(x,y)&=x^py^p\log\big(\tfrac{y}{x}\big)
   \qquad&&\mbox{if}\quad p=q\in\R,\\
  \D_{f,g}(x,y)&=x^ay^a\sin\big(b\log\big(\tfrac{x}{y}\big)\big) \qquad&&\mbox{if}\quad
      p=\bar{q}=a+bi,\,\,a,b\in\R,\,b\neq0.
}
In the first two cases, $\D_{f,g}(x,y)$ is different from zero for all $x,y\in \R_+$ with $x\neq y$, 
therefore $f$ and $g$ form a Chebyshev system on $\R_+$ if $p,q\in\R$.  In the third case, that is if 
$p=\bar{q}\not\in\R$, the functions $f$ and $g$ may form a Chebyshev system only in a subinterval $I$ of 
$\R_+$ only. One can easily see that then $\D_{f,g}(x,y)$ is different from zero for all $x,y\in I$ with 
$x\neq y$ if and only if the interval $I/I:=\{u/v\mid u,v\in I\}$ is contained in the open interval 
$\big(\exp\big(-\tfrac{\pi}{|b|}\big),\exp\big(\tfrac{\pi}{|b|}\big)\big)$. Obviously this cannot happen if 
either $\inf I=0$ or $\sup I=+\infty$. The $d$-variable generalized Gini mean $G_{p,q,m;\mu}$ is defined now 
to be the $d$-variable Bajraktarevi\'c mean $M_{f,g,m;\mu}$, where $f$ and $g$ are given by \eq{pq}.

Let $I$ be an open subinterval of $\R_+$ such that $I$ is contained in 
$\big(\exp\big(-\tfrac{\pi}{2|b|}\big),\exp\big(\tfrac{\pi}{2|b|}\big)\big)$ if $p=\bar{q}\not\in\R$. Then 
$g$ is positive over $I$ in each of the three possibilities and for $\pmb{x}\in I^d$, and, as a particular 
case of equation \eq{fg}, we can obtain the following explicit formula for the $d$-variable Gini mean 
$G_{p,q,m;\mu}$:
\Eq{ab}{
   G_{p,q,m;\mu}(\pmb{x})
      :=\begin{cases}
          \left(\dfrac{\int_T \big(m(\pmb{x},t)\big)^p\d\mu(t)}
                 {\int_T \big(m(\pmb{x},t)\big)^q\d\mu(t)}\right)^{\frac1{p-q}} 
                 &\mbox{if }p,q\in\R,\,p\neq q,\\[5mm]
          \exp\left(\dfrac{\int_T \big(m(\pmb{x},t)\big)^p\log\big(m(\pmb{x},t)\big)\d\mu(t)}
                 {\int_T \big(m(\pmb{x},t)\big)^p\d\mu(t)}\right) &\mbox{if }p=q\in\R,\\[5mm]
          \exp\left(\dfrac{1}{b}\arctan\left(\dfrac{\int_T 
                 \big(m(\pmb{x},t)\big)^a\sin\big(b\log\big(m(\pmb{x},t)\big)\big)\d\mu(t)}          
         {\int_T\big(m(\pmb{x},t)\big)^a\cos\big(b\log\big(m(\pmb{x},t)\big)\big)\d\mu(t)}\right)\right) 
                 &\mbox{if }p=\bar{q}=a+bi\not\in\R.
        \end{cases}
}
If $q=0$ in the above formula, we get the definition of $d$-variable \textit{generalized H\"older means} as follows:
\Eq{a}{
   G_{p,0,m;\mu}(\pmb{x}):=H_{p,m;\mu}(\pmb{x})
      :=\begin{cases}
          \left(\dfrac{\int_T \big(m(\pmb{x},t)\big)^p\d\mu(t)}
                 {\int_T \d\mu(t)}\right)^{\frac1{p}} &\mbox{if }p\in\R\setminus\{0\},\\[5mm]
          \exp\left(\dfrac{\int_T \log\big(m(\pmb{x},t)\big)\d\mu(t)}
                 {\int_T \d\mu(t)}\right) &\mbox{if }p=0,
        \end{cases}
   \qquad(\pmb{x}\in \R_+^d).
}

In the particular case when $T=\{1,\dots,d\}$, $\mu=\frac{\delta_1+\dots+\delta_d}{d}$ and
$m(\pmb{x},t)=x_t$, formula \eq{ab} reduces to the so-called $d$-variable \emph{Gini mean} $G_{p,q}$
(cf.\ \cite{Gin38}): for $\pmb{x}=(x_1\dots,x_d)\in I^d$,
\Eq{*}{
   G_{p,q}(\pmb{x})
      :=\begin{cases}
         \bigg(\dfrac{x_1^p+\dots+x_d^p}{x_1^q+\dots+x_d^q}\bigg)^{\frac1{p-q}}
         &\mbox{if }p,q\in\R,\,p\neq q,\\[5mm]
         \exp\bigg(\dfrac{x_1^p\log(x_1)+\dots+x_d^p\log(x_d)}{x_1^p+\dots+x_d^p}\bigg)
         &\mbox{if }p=q\in\R, \\[5mm]
         \exp\left(\dfrac{1}{b}\arctan\left(\dfrac{x_1^a\sin(\log(x_1^b))+\cdots+x_d^a\sin(\log(x_d^b))}
        {x_1^a\cos(\log(x_1^b))+\cdots+x_d^a\cos(\log(x_d^b))}\right)\right) 
          &\mbox{if }p=\bar{q}=a+bi\not\in\R.
        \end{cases}
}
Gini means in this generality (i.e., including the case of non-real parameters $p,q$ were dealt with in the 
paper \cite{Pal89c} where the comparison problem of these means was solved.
Obviously, $G_{p,0}=H_p$, i.e., H\"older means are particular Gini means. For further particular cases of
formula \eq{fg}, we refer to the paper \cite{PalZak17}.

The aim of this paper is to study the equality and the homogeneity problems of these means, i.e., to find
conditions for the generating functions $(f,g)$ and $(h,k)$, for the family of means $m$, and for
the measure $\mu$ such that the functional equation
\Eq{2}{
   M_{f,g,m;\mu}(\pmb{x})=M_{h,k,m;\mu}(\pmb{x}) \qquad(\pmb{x}\in I^d)
} 
and the homogeneity property
\Eq{*}{
   M_{f,g,m;\mu}(\lambda\pmb{x})=\lambda M_{f,g,m;\mu}(\pmb{x}) 
   \qquad(\lambda>0,\,\pmb{x},\lambda\pmb{x}\in I^d),
}
respectively, be satisfied. Our main results generalize that of the paper by Losonczi and P\'ales
\cite{LosPal11a}, Losonczi \cite{Los13} and also many former results obtained in various particular cases of 
this problem, cf.\ \cite{AczDar63c}, \cite{BerMor98}, \cite{BerMor00}, \cite{BurJar13}, \cite{Los99}, 
\cite{Los02d}, \cite{Los02a}, \cite{Los02c}, \cite{Los03a}, \cite{Los06b}, \cite{Los07a}, \cite{Los07b},
\cite{MakPal08}, \cite{Pal11}. As direct applications of the results obtained on the equality of generalized 
Bajraktarevi\'c means, we consider and solve the homogeneity problem of these means under general conditions.

\section{Auxiliary results}
\setcounter{equation}{0} 

In order to describe the regularity conditions related to the two generating functions $f,g$ of the mean
$M_{f,g,m;\mu}$, we introduce some regularity classes. The class $\C_0(I)$ consists of all those pairs of
continuous functions $f,g:I\to\R$ that form a Chebyshev system over $I$.

If $n\ge1$, then we say that the pair $(f,g)$ is in the class $\C_n(I)$
if $f,g$ are $n$-times continuously differentiable functions such that $(f,g)\in\C_0(I)$ and the
\textit{Wronski determinant}
\Eq{W}{
  \DET{f'(x)&f(x)\\g'(x)&g(x)}
     =\partial_1\D_{f,g}(x,x)
   \qquad(x\in I)
}
does not vanish on $I$. Provided that $g$ is positive, then we have that
\Eq{f/g}{
  \bigg(\frac{f}{g}\bigg)'(x)=\frac{\partial_1\D_{f,g}(x,x)}{g^2(x)}
}
hence condition $\partial_1\D_{f,g}(x,x)\neq0$ implies that $f/g$ is
strictly monotone, whence it follows that $(f,g)\in\C_0(I)$.
Obviously, $\C_0(I)\supseteq\C_1(I)\supseteq\C_2(I)\supseteq\cdots$.

It is easy to see that if $(f,g),(h,k)\in\C_0(I)$ and
\Eq{5}{
  f&=\alpha h+\beta k,\\
  g&=\gamma h+\delta k,
}
where the constants $\alpha, \beta, \gamma, \delta\in \R$ satisfy
$\alpha\delta-\beta\gamma\ne 0$, then, by the product theorem for determinants, it follows that
\Eq{4.5}{
  \D_{f,g}=\DET{\alpha&\beta\\\gamma&\delta}\cdot\D_{h,k}.
}
This, in view of \lem{0}, implies that the identity
\Eq{4}{
  M_{f,g,m;\mu}=M_{h,k,m;\mu}
}
also holds for any measurable family of $d$-variable means $m:I^d\times T\to I$ and probability measure $\mu$.

If \eq{5} holds for some constants $\alpha, \beta, \gamma, \delta\in \R$ with $\alpha\delta-\beta\gamma\ne 0$,
then we say that the pairs $(f,g)$ and $(h,k)$ are \emph{equivalent}. It is obvious that any necessary and/or
sufficient condition for \eq{2} has to be invariant with respect to the equivalence of the generating
functions.

For the characterization of the equivalence, we introduce the following notations:
for $(f,g)\in\C_2(I)$, the functions $\Phi_{f,g},\Psi_{f,g}:I\to\R$ are defined by
\Eq{D1}{
\Phi_{f,g}(x):=\frac{\partial_1^2\D_{f,g}(x,x)}{\partial_1\D_{f,g}(x,x)}
\qquad\mbox{and}\qquad
\Psi_{f,g}(x):=-\frac{\partial_1^2\partial_2\D_{f,g}(x,x)}{\partial_1\D_{f,g}(x,x)}\qquad(x\in I).
}
In other words,
\Eq{DD}{
  \Phi_{f,g}:=\frac{\DET{f''&f\\g''&g}}{\DET{f'&f\\g'&g}}
  \qquad\mbox{and}\qquad
  \Psi_{f,g}:=-\frac{\DET{f''&f'\\g''&g'}}{\DET{f'&f\\g'&g}}.
}

\Thm{1.5}{If $(f,g),(h,k)\in\C_2(I)$, then the pairs $(f,g)$ and $(h,k)$ are equivalent if and only if
\Eq{1.5}{
  \Phi_{f,g}=\Phi_{h,k}\qquad\mbox{and}\qquad \Psi_{f,g}=\Psi_{h,k}.
}}

\begin{proof}
If $(f,g)$ and $(h,k)$ are equivalent, then, for some $\alpha, \beta, \gamma, \delta\in \R$ with
$\alpha\delta-\beta\gamma\ne 0$, we have \eq{5}, which implies \eq{4.5}. Using this formula and the
definition of $\Phi_{f,g},\Psi_{f,g},\Phi_{h,k}$, and $\Psi_{h,k}$, the identities in \eq{1.5}
follow directly.

Now assume that \eq{1.5} is valid on $I$ and consider the following second order homogeneous linear
differential equation
\Eq{y''}{
   y''=\Phi_{f,g}y'+\Psi_{f,g}y.
}
Using the definitions of $\Phi_{f,g},\Psi_{f,g}$ from \eq{DD}, we can rewrite it in the following equivalent
form
\Eq{*}{
  y''=\frac{\DET{f''&f\\g''&g}}{\DET{f'&f\\g'&g}}y'-\frac{\DET{f''&f'\\g''&g'}}{\DET{f'&f\\g'&g}}y.
}
After multiplying this equation by $\DET{f'&f\\g'&g}$, and rearranging every term to one side of the
equation, we infer that \eq{y''} is equivalent to
\Eq{EE}{
  \DET{y''&y'&y\\f''&f'&f\\g''&g'&g}=0.
}
The functions $y=f$ and $y=g$ are trivially solutions of \eq{EE}, therefore they are solutions of \eq{y''} as
well. Their Wronski determinant is nonzero, hence every solution of \eq{y''} is a linear combination of them.

On the other hand, due to the identities \eq{1.5}, the differential equation \eq{y''} is also equivalent to
\Eq{*}{
   y''=\Phi_{h,k}y'+\Psi_{h,k}y.
}
By a similar argument as above, we can see that $h$ and $k$ are also linearly independent solutions of this
second order homogeneous linear differential equation. Therefore, $f$ and $g$ should be their
(independent) linear combinations, i.e., \eq{5} should hold for some $\alpha, \beta, \gamma, \delta\in \R$
with $\alpha\delta-\beta\gamma\ne 0$. This proves the equivalence of the pairs $(f,g)$ and $(h,k)$.
\end{proof}

The following lemma is an immediate consequence of the asymmetry property
\Eq{*}{
  \D_{f,g}(x,y)=-\D_{f,g}(y,x) \qquad(x,y\in I).
}

\Lem{1.5}{If $(f,g)\in\C_n(I)$ for $n\in\{1,2,3\}$ then, for all $x\in I$,
\Eq{id12}{
\partial_2\D_{f,g}(x,x)&=-\partial_1\D_{f,g}(x,x),\qquad &
\partial_1\partial_2\D_{f,g}(x,x)&=0, \\
\partial_2^2\D_{f,g}(x,x)&=-\partial_1^2\D_{f,g}(x,x), \qquad &
\partial_2^2\partial_1\D_{f,g}(x,x)&=-\partial_1^2\partial_2\D_{f,g}(x,x),
}
and
\Eq{id3}{
\partial_2^3\D_{f,g}(x,x)&=-\partial_1^3\D_{f,g}(x,x),
}
respectively.}

\Lem{2-}{If $(f,g)\in\C_3(I)$, then
\Eq{D2}{
\frac{\partial_1^3\D_{f,g}(x,x)}{\partial_1\D_{f,g}(x,x)}=\Phi'_{f,g}(x)+\Phi_{f,g}^2(x)+\Psi_{f,g}
(x)\qquad(x\in I).
}}

\begin{proof}
By computing the derivative of $\Phi_{f,g}(x)$, we get
\Eq{*}{
\Phi'_{f,g}(x)=\frac{\partial_1^3\D_{f,g}(x,x)+\partial_2\partial_1^2\D_{f,g}(x,x)}{\partial_1\D_{f,g}(x,x)}
-\frac{\partial_1^2\D_{f,g}(x,x)\Big(\partial_1^2\D_{f,g}(x,x)+\partial_2\partial_1\D_{f,g}(x,x)\Big)}{
(\partial_1\D_{f,g}(x,x))^2}.
}
Since $\partial_1\partial_2\D_{f,g}(x,x)=0$ is consequence of the asymmetry property
$\D_{f,g}(x,y)=-\D_{f,g}(y,x)$, we get that
\Eq{*}{
\Phi'_{f,g}(x)
  =\frac{\partial_1^3\D_{f,g}(x,x)}{\partial_1\D_{f,g}(x,x)}
     +\frac{\partial_2\partial_1^2\D_{f,g}(x,x)}{\partial_1\D_{f,g}(x,x)}
     -\bigg(\frac{\partial_1^2\D_{f,g}(x,x)}{\partial_1\D_{f,g}(x,x)}\bigg)^2
  =\frac{\partial_1^3\D_{f,g}(x,x)}{\partial_1\D_{f,g}(x,x)}
     -\Psi_{f,g}(x)-\Phi^2_{f,g}(x),
}
whence equation \eq{D2} follows immediately.
\end{proof}

The following result, which is based on \cite[Theorem 3]{BesPal03}, allows us to assume more regularity on
Chebyshev systems.

\Lem{BP}{Let $n\in\N\cup\{0\}$ and $(f,g)\in\C_n(I)$. Then there exist $\alpha, \beta, \gamma, \delta\in \R$
with $\alpha\delta-\beta\gamma\ne 0$ and $(h,k)\in\C_n(I)$ such that \eq{5} holds and $k$ is positive
and $h/k$ is strictly monotone. Furthermore, if $n\geq1$, then the derivative of $h/k$ does not
vanish on $I$.}

For its proof, the reader should consult \cite[Theorem 3]{BesPal03} and \cite[Lemma 2]{PalZak17} for $n=0$
and $n\geq1$, respectively.

For the computation of the first-, second- and third-order partial derivatives of the mean $M_{f,g,m;\mu}$ at
the
diagonal of $I^d$, we will establish a result below. For brevity, we introduce the following notation:
If $\pmb{p}\in I^d$ and $\delta>0$ then let $B(\pmb{p},\delta)$ stand for the ball $\{\pmb{x}\in I^d\colon
|\pmb{x}-\pmb{p}|\leq\delta\}$. Furthermore, if $\mu$ is a probability measure on the measurable
space $(T,\A)$ and $q\geq1$, then the space of measurable functions $\varphi:T\to\R$ such that
$|\varphi|^q$ is $\mu$-integrable will be denoted by $L^q(T,\A,\mu)$ or shortly by $L^q$.

If $\varphi:T\to\R$ is a $\mu$-integrable function, then we set
\Eq{*}{
  \langle\varphi\rangle_\mu:=\int_T\varphi(t)\d\mu(t).
}
More generally, if $\varphi:I^d\times T\to\R$, and for some $\pmb{x}\in I^d$, the map
$t\mapsto\varphi(\pmb{x},t),t)$ is $\mu$-integrable, then we write
\Eq{PD+}{
  \langle\varphi\rangle_\mu(\pmb{x}):=\int_T\varphi(\pmb{x},t)\d\mu(t).
}
Given a number $q\geq1$, a function $\varphi:I^d\times T\to\R$ is said to be of $L^q$-type at $\pmb{p}\in
I^d$, if $\varphi(\pmb{p},\cdot)$ is measurable, furthermore, there exist $\delta>0$ and a function $a\in L^q$
such that
\Eq{*}{
  |\varphi(\pmb{x},t)|\leq a(t) \qquad(t\in T,\,\pmb{x}\in B(\pmb{p},\delta)).
}

Let $\C_1(I^d\times T)$ denote the class of measurable families of $d$-variable means $m:I^d\times T\to I$
with the following two additional properties:
\begin{enumerate}
 \item[(H5)] For every $t\in T$, the function $m(\cdot,t)$ is continuously partially differentiable over
$I^d$ such that, for all $\pmb{p}\in I^d$, $i\in\{1,\dots,d\}$, the function $\partial_i m$ is of
$L^1$-type at $\pmb{p}$.
\end{enumerate}
Analogously, we define $\C_2(I^d\times T)$ to be the following subclass of $\C_1(I^d\times T)$:
\begin{enumerate}
 \item[(H6)] For every $t\in T$, the function $m(\cdot,t)$ is twice continuously partially differentiable over
$I^d$ such that, for all $\pmb{p}\in I^d$ and $i,j\in\{1,\dots,d\}$, the function $\partial_i m$ is of
$L^2$-type and $\partial_i \partial_j m$ is of $L^1$-type at $\pmb{p}$.
\end{enumerate}
Similarly, we define $\C_3(I^d\times T)$ to be the following subclass of $\C_2(I^d\times T)$:
\begin{enumerate}
\item[(H7)]For every $t\in T$, the function $m(\cdot,t)$ is three times continuously partially differentiable over
$I^d$ such that, for all $\pmb{p}\in I^d$ and $i,j,l\in\{1,\dots,d\}$, the function $\partial_i m$ is of
$L^3$-type, $\partial_i \partial_j m$ is of $L^{\frac{3}{2}}$-type, and $\partial_i \partial_j \partial_l m$ is of $L^1$-type at $\pmb{p}$.
\end{enumerate}

In order to formulate the results below, we introduce the following notation: for $i,j,l\in\{1,\dots,d\}$,
define $\sigma(i,j,l)$ to be the set of all cyclic permutations of $(i,j,l)$, that is,
\Eq{*}{
  \sigma(i,j,l):=\{(i,j,l),(j,l,i),(l,i,j)\}.
}

\Lem{1-}{Let $n\in\{1,2,3\}$ and let $\varphi:I\to\R$ be a $n$-times continuously differentiable function and
$m\in \C_n(I^d\times T)$. Then the function $\Phi:I^d\to\R$ defined by
\Eq{1.0}{
 \Phi(\pmb{x}):=\int_T\varphi(m(\pmb{x},t))\d\mu(t)
}
is $n$-times continuously differentiable on $I^d$. Furthermore, for $i\in\{1,\dots,d\}$,
\Eq{1.1}{
  \partial_i \Phi(\pmb{p})=\int_T\varphi'(m(\pmb{p},t))\,\partial_im(\pmb{p},t)\d\mu(t) \qquad(\pmb{p}\in I^d)
}
for $i,j\in\{1,\dots,d\}$ and $n=2$,
\Eq{1.2}{
  \partial_i \partial_j \Phi(\pmb{p})
   =\int_T\Big[\varphi''(m(\pmb{p},t))\,\partial_im(\pmb{p},t)\,\partial_jm(\pmb{p},t)
   +\varphi'(m(\pmb{p},t))\,\partial_i\partial_jm(\pmb{p},t)\Big]\d\mu(t) \qquad(\pmb{p}\in I^d)
}
and, for $i,j,l\in\{1,\dots,d\}$ and $n=3$,
\Eq{1.3}{
  \partial_i \partial_j\partial_l \Phi(\pmb{p})
=\int_T\Big[\varphi'''(m(\pmb{p},t))\prod_{r\in\{i,j,l\}}\partial_rm(\pmb{p},t)
&+\varphi''(m(\pmb{p},t))\!\!\!\sum_{(\alpha,\beta,\gamma)\in\sigma(i,j,l)}\!\!\!
\big(\partial_\alpha m(\pmb{p},t)\,\partial_\beta\partial_\gamma m(\pmb{p},t)\\
     &+\varphi'(m(\pmb{p},t))\,\partial_i\partial_j\partial_lm(\pmb{p},t)\Big]\d\mu(t)
\qquad(\pmb{p}\in I^d).
}
}

The proof of the above lemma for the cases $n=1$ and $n=2$ was elaborated in details in the paper
\cite{PalZak17}, the argument concerning the case $n=3$ is completely analogous, therefore, we omit it.
For the sake of convenience, introduce the following notations: For $m\in \C_1(I^d\times T)$ and
$r\in\{1,\dots,d\}$, denote
\Eq{*}{
  \partial_r^*m(\pmb{x},t):=\partial_r m(\pmb{x},t)-\langle \partial_r m\rangle_\mu (\pmb{x})
  \qquad(\pmb{x}\in I^d,\,t\in T),
}
and, for $x\in I$, set
\Eq{*}{
  x^{(d)}:=(x,\dots,x)\in I^d.
}

\Thm{1+}{Let $(f,g)\in\C_1(I)$, let $m\in\C_1(I^d\times T)$ be a measurable family of means, and let $\mu$ be
a probability measure on the measurable space $(T,\A)$.
Then $M_{f,g,m;\mu}$ is continuously differentiable on $I^d$ and, for all $i\in\{1,\dots,d\}$ and $x\in I$,
\Eq{PD1}{
  \partial_iM_{f,g,m;\mu}\big(x^{(d)}\big)
  =\langle\partial_im\rangle_\mu\big(x^{(d)}\big).
}
If, in addition, $(f,g)\in\C_2(I)$, let $m\in\C_2(I^d\times T)$, then $M_{f,g,m;\mu}$ is twice continuously
differentiable on $I^d$ and, for all $i,j\in\{1,\dots,d\}$ and $x\in I$,
\Eq{PD2}{
  \partial_i\partial_j &M_{f,g,m;\mu}\big(x^{(d)}\big)
  =\Phi_{f,g}(x)\big\langle \partial_i^*m \,\partial_j^*m\big\rangle_\mu\big(x^{(d)}\big)
    +\langle \partial_i\partial_j m\rangle_\mu\big(x^{(d)}\big).
}
Finally, if  $(f,g)\in\C_3(I)$, let $m\in\C_3(I^d\times T)$, then $M_{f,g,m;\mu}$ is three times continuously
differentiable on $I^d$ and, for all $i,j,l\in\{1,\dots,d\}$ and $x\in I$,
\Eq{PD3}{
  \partial_i\partial_j\partial_l M_{f,g,m;\mu}\big(x^{(d)}\big)
  =&\big(\Phi'_{f,g}(x)+\Phi_{f,g}^2(x)\big)
  \big(\big\langle\partial_im\,\partial_jm\,\partial_lm\big\rangle_\mu
-\langle\partial_im\rangle_\mu\langle\partial_jm\rangle_\mu\langle\partial_lm\rangle_\mu\big)\big(x^{(d)}
\big)\\
  &+\Phi_{f,g}(x)\sum_{(\alpha,\beta,\gamma)\in\sigma(i,j,l)}\!\!\!
   \big(\langle\partial_\alpha \partial_\beta m\,\partial_\gamma m\rangle_\mu
   -\partial_\alpha\partial_\beta M_{f,g,m;\mu}\partial_\gamma M_{f,g,m;\mu}\big)\big(x^{(d)}\big)\\
&+\Psi_{f,g}(x)\big\langle\partial_i^*m\,\partial_j^*m\,\partial_l^*m\big\rangle_\mu\big(x^{(d)}\big)
  + \langle\partial_i\partial_j\partial_lm\rangle_\mu\big(x^{(d)}\big).
}}
\begin{proof} Let $n\in\{1,2,3\}$ and assume that $(f,g)\in\C_n(I)$, $m\in\C_n(I^d\times T)$. In view of
\lem{BP}, we may assume that $g$ is positive, $f/g$ is strictly monotone with a non-vanishing first-order
derivative. Then $f$, $g$ and the inverse of $f/g$ are $n$-times continuously differentiable and, by \lem{1-},
we also have that the mappings
\Eq{*}{
  \pmb{x}\mapsto \int_T f\big(m(\pmb{x},t)\big)\d\mu(t)
  \qquad\mbox{and}\qquad
  \pmb{x}\mapsto \int_T g\big(m(\pmb{x},t)\big)\d\mu(t)
}
are $n$-times continuously differentiable on $I^d$. On the other hand, we now also
have formula \eq{fg} for the $d$-variable mean $M_{f,g,m;\mu}$. Thus, using the standard calculus rules, it
follows that $M_{f,g,m;\mu}$ is $n$-times continuously differentiable on $I^d$.

For the equalities \eq{PD1} and \eq{PD2}, the reader shall establish \cite[Theorem 4]{PalZak17}. To prove the
third formula stated in \eq{PD3}, let us consider the case $n=3$. In view of \lem{0}, we have the following
identity
\Eq{DM}{
  \int_T \D_{f,g}(m(\pmb{x},t),M_{f,g,m;\mu}(\pmb{x}))\d\mu(t)=0 \qquad (\pmb{x}\in I^d).
}
For the sake of brevity, for $\alpha,\beta\in\N\cup\{0\}$ with $1\leq\alpha+\beta\leq3$, introduce the
following notation
\Eq{*}{
  \Delta_{\alpha,\beta}(\pmb{x},t)
  :=\partial_1^\alpha\partial_2^\beta\D_{f,g}(m(\pmb{x},t),M_{f,g,m;\mu}(\pmb{x}))
  \qquad (\pmb{x}\in I^d\,t\in T).
}
Performing the partial differentiations $\partial_i$, $\partial_j\partial_i$ and
$\partial_l\partial_j\partial_i$ on equality \eq{DM} side by side, we get
\Eq{*}{
   \int_T \big[\Delta_{1,0}(\pmb{x},t)\partial_im(\pmb{x},t)
+\Delta_{0,1}(\pmb{x},t)\partial_iM_{f,g,m;\mu}(\pmb{x})\big]\d\mu(t)=0,
}
\Eq{*}{
   \int_T \big[\Delta_{2,0}(\pmb{x},t)
             &\partial_im(\pmb{x},t)\partial_jm(\pmb{x},t)
   +\Delta_{0,2}(\pmb{x},t)
             \partial_jM_{f,g ,m;\mu}(\pmb{x})\partial_iM_{f,g ,m;\mu}(\pmb{x})\\
  &+\Delta_{1,1}(\pmb{x},t)
             \big(\partial_jM_{f,g,m;\mu}(\pmb{x})\partial_im(\pmb{x},t)
             +\partial_iM_{f,g,m;\mu}(\pmb{x})\partial_jm(\pmb{x},t)\big)\\
  &+\Delta_{1,0}(\pmb{x},t)\partial_j\partial_im(\pmb{x},t)
  +\Delta_{0,1}(\pmb{x},t)\partial_j\partial_iM_{f,g,m;\mu}(\pmb{x})\big]\d\mu(t)=0,
}
and
\Eq{*}{
  \int_T &\Big[\Delta_{3,0}(\pmb{x},t)\prod_{r\in\{i,j,l\}}\partial_rm(\pmb{x},t)
  +\Delta_{0,3}(\pmb{x},t)\prod_{r\in\{i,j,l\}}\partial_r M_{f,g,m;\mu}(\pmb{x})\\
&+\Delta_{2,1}(\pmb{x},t)\!\!\!\sum_{(\alpha,\beta,\gamma)\in\sigma(i,j,l)}\!\!\!
   \partial_\alpha M_{f,g,m;\mu}(\pmb{x})\partial_\beta m(\pmb{x},t)\partial_\gamma m(\pmb{x},t)\\
&+\Delta_{1,2}(\pmb{x},t)\!\!\!\sum_{(\alpha,\beta,\gamma)\in\sigma(i,j,l)}\!\!\!
   \partial_\alpha M_{f,g,m;\mu}(\pmb{x})\partial_\beta M_{f,g,m;\mu}(\pmb{x})\partial_\gamma m(\pmb{x},t)\\
&+\Delta_{2,0}(\pmb{x},t)\!\!\!\sum_{(\alpha,\beta,\gamma)\in\sigma(i,j,l)}\!\!\!
   \partial_\alpha \partial_\beta m(\pmb{x},t)\partial_\gamma m(\pmb{x},t)
   +\Delta_{0,2}(\pmb{x},t)\!\!\!\sum_{(\alpha,\beta,\gamma)\in\sigma(i,j,l)}\!\!\!
   \partial_\alpha\partial_\beta M_{f,g,m;\mu}(\pmb{x})\partial_\gamma M_{f,g,m;\mu}(\pmb{x})\\
&+\Delta_{1,1}(\pmb{x},t)\!\!\!\sum_{(\alpha,\beta,\gamma)\in\sigma(i,j,l)}\!\!\!
   \Big(\partial_\alpha\partial_\beta M_{f,g,m;\mu}(\pmb{x})\partial_\gamma m(\pmb{x},t)
    +\partial_\alpha M_{f,g,m;\mu}(\pmb{x})\partial_\beta \partial_\gamma m(\pmb{x},t)\Big)\\
&+\Delta_{1,0}(\pmb{x},t)\partial_i\partial_j\partial_lm(\pmb{x},t)
+\Delta_{0,1}(\pmb{x},t)\partial_i\partial_j\partial_lM_{f,g,m;\mu}(\pmb{x})\Big]\d\mu(t)=0,
}
respectively. Using the identities \eq{id12}, \eq{id3} and substituting $x^{(d)}\in I^d$, we get that
\Eq{*}{
  \int_T &\Big[\partial_1^3\D_{f,g}(x,x)\Big(\prod_{r\in\{i,j,l\}}\partial_rm\big(x^{(d)},t\big)
  -\prod_{r\in\{i,j,l\}}\partial_r M_{f,g,m;\mu}\big(x^{(d)}\big)\Big)\\
&+\partial_1^2\partial_2\D_{f,g}(x,x)\!\!\!\sum_{(\alpha,\beta,\gamma)\in\sigma(i,j,l)}\!\!\!
   \partial_\alpha M_{f,g,m;\mu}\big(x^{(d)}\big)\Big(\partial_\beta m\big(x^{(d)},t\big)
 -\partial_\beta M_{f,g,m;\mu}\big(x^{(d)}\big)\Big)\partial_\gamma m\big(x^{(d)},t\big)\\
&+\partial_1^2\D_{f,g}(x,x)\!\!\!\sum_{(\alpha,\beta,\gamma)\in\sigma(i,j,l)}\!\!\!
   \Big(\partial_\alpha \partial_\beta m\big(x^{(d)},t\big)\partial_\gamma m\big(x^{(d)},t\big)
   -\partial_\alpha\partial_\beta M_{f,g,m;\mu}\big(x^{(d)}\big)\partial_\gamma
M_{f,g,m;\mu}\big(x^{(d)}\big)\Big)\\
&+\partial_1\D_{f,g}(x,x)\partial_i\partial_j\partial_lm\big(x^{(d)},t\big)
-\partial_1\D_{f,g}(x,x)\partial_i\partial_j\partial_lM_{f,g,m;\mu}\big(x^{(d)}\big)\Big]\d\mu(t)=0.
}
Now, dividing the above equality by $\partial_1\D_{f,g}(x,x)$ and using the definitions of $\Phi_{f,g}$ and
$\Psi_{f,g}$, and also the identity \eq{D2}, for $x^{(d)}\in I^d$, the following formula follows
\Eq{*}{
  \partial_i&\partial_j\partial_l  M_{f,g,m;\mu}\big(x^{(d)}\big) \\
  =&\big(\Phi'_{f,g}(x)+\Phi_{f,g}^2(x)+\Psi_{f,g}(x)\big)
  \int_T \Big(\prod_{r\in\{i,j,l\}}\partial_rm\big(x^{(d)},t\big)
  -\prod_{r\in\{i,j,l\}}\partial_r M_{f,g,m;\mu}\big(x^{(d)}\big)\Big)\d\mu(t)\\
  &-\Psi_{f,g}(x)\!\!\!\sum_{(\alpha,\beta,\gamma)\in\sigma(i,j,l)}\!\!\!\partial_\alpha
   M_{f,g,m;\mu}\big(x^{(d)}\big)\int_T \Big(\partial_\beta m\big(x^{(d)},t\big)
   -\partial_\beta M_{f,g,m;\mu}\big(x^{(d)}\big)\Big)\partial_\gamma m\big(x^{(d)},t\big)\d\mu(t) \\
  &+\Phi_{f,g}(x)\!\!\!\sum_{(\alpha,\beta,\gamma)\in\sigma(i,j,l)}\int_T
   \Big(\partial_\alpha \partial_\beta m\big(x^{(d)},t\big)\partial_\gamma m\big(x^{(d)},t\big)
   -\partial_\alpha\partial_\beta M_{f,g,m;\mu}\big(x^{(d)}\big)\partial_\gamma
M_{f,g,m;\mu}\big(x^{(d)}\big)\Big)\d\mu(t)\\
  &+ \int_T\partial_i\partial_j\partial_lm\big(x^{(d)},t\big)\d\mu(t).
}
By using formula \eq{PD1} and definition \eq{PD+}, we get the following simplified equation
\Eq{*}{
  \partial_i\partial_j\partial_l & M_{f,g,m;\mu}\big(x^{(d)}\big) \\
  =&\big(\Phi'_{f,g}(x)+\Phi_{f,g}^2(x)\big)
  \big(\big\langle\partial_im\,\partial_jm\,\partial_lm\big\rangle_\mu
-\langle\partial_im\rangle_\mu\langle\partial_jm\rangle_\mu\langle\partial_lm\rangle_\mu\big)\big(x^{(d)}
\big)\\
  &+\Phi_{f,g}(x)\sum_{(\alpha,\beta,\gamma)\in\sigma(i,j,l)}\!\!\!
   \Big(\langle\partial_\alpha \partial_\beta m\,\partial_\gamma m\rangle_\mu
   -\partial_\alpha\partial_\beta M_{f,g,m;\mu}\partial_\gamma M_{f,g,m;\mu}\Big)\big(x^{(d)}\big)\\
   &+\Psi_{f,g}(x)\Big(\big\langle\partial_im\,\partial_jm\,\partial_lm\big\rangle_\mu
  -\!\!\!\sum_{(\alpha,\beta,\gamma)\in\sigma(i,j,l)}\!\!\!\langle\partial_\alpha m\rangle_\mu
\langle\partial_\beta m\,\partial_\gamma m\rangle_\mu
+2\langle\partial_im\rangle_\mu\langle\partial_jm\rangle_\mu\langle\partial_lm\rangle_\mu\Big)\big(x^{(d)}
\big)\\
  &+ \langle\partial_i\partial_j\partial_lm\rangle_\mu\big(x^{(d)}\big).
}
This equality, combined with the following easy-to-see identity,
\Eq{*}{
\big\langle\partial_im\,\partial_jm\,\partial_lm\big\rangle_\mu
  -\!\!\!\sum_{(\alpha,\beta,\gamma)\in\sigma(i,j,l)}\!\!\!\langle\partial_\alpha m
\rangle_\mu\langle\partial_\beta m\,\partial_\gamma m\rangle_\mu
+2\langle\partial_im\rangle_\mu\langle\partial_jm\rangle_\mu\langle\partial_lm\rangle_\mu
=\big\langle\partial_i^*m\,\partial_j^*m\,\partial_l^*m\big\rangle_\mu
}
yields formula \eq{PD3} of the theorem.
\end{proof}

We note that, using the equality in \eq{PD2}, the formula \eq{PD3} for the third-order partial derivatives
can be made more explicit.

\section{Equality of generalized Bajraktarevi\'c means}\setcounter{equation}{0} 

In this section we characterize the equality of generalized Bajraktarevi\'c and quasi-arithmetic means under
3 times and 2 times differentiability assumptions, respectively.

\Thm{MT}{Let $(f,g),(h,k)\in\C_3(I)$, let $m\in\C_3(I^d\times T)$ be a measurable family of means, and let
$\mu$ be a probability measure on the measurable space $(T,\A)$. Assume that, there exists a dense subset
$D\subseteq I$ such that, for all $x\in D$,
\Eq{MT-1}{
  \mu\Big(\big\{t\in T\mid
\partial_1^*m\big(x^{(d)},t\big)=\cdots=\partial_d^*m\big(x^{(d)},t\big)=0\big\}\Big)<1
}
and there exist $i,j,l\in\{1,\dots,d\}$ such that
\Eq{MT0}{
\big\langle\partial_i^*m\,\partial_j^*m\,\partial_l^*m\big\rangle_\mu\big(x^{(d)}\big)\neq0.
}
Then the following assertions are equivalent:
\begin{enumerate}[(i)]
 \item For all $\pmb{x}\in I^d$,
 \Eq{MT1}{
   M_{f,g,m;\mu}(\pmb{x})=M_{h,k,m;\mu}(\pmb{x});
 }
 \item There exists an open set $U\subseteq I^d$ containing the subdiagonal $\{x^{(d)}\mid x\in D\}$
such that, for all $\pmb{x}\in U$, the equality \eq{MT1} holds;
 \item The two identities in \eq{1.5} hold;
 \item The pairs $(f,g)$ and $(h,k)$ are equivalent.
\end{enumerate}}

\begin{proof} The implication \textit{(i)}$\Rightarrow$\textit{(ii)} is trivial. As we have seen it at the
beginning of Section 2, equivalent pairs generate identical means, hence the implication
\textit{(iv)}$\Rightarrow$\textit{(i)} is also valid. The implication
\textit{(iii)}$\Rightarrow$\textit{(iv)} is the consequence of \thm{1.5}. Therefore, it remains to show that
condition \textit{(ii)} implies \textit{(iii)}.

Assume that condition \textit{(ii)} holds for some open set $U\subseteq I^d$ containing the subdiagonal
$\{x^{(d)}\mid x\in D\}$. By the regularity assumptions of the theorem, the two means are 3 times
continuously differentiable over $U$. Therefore, for all $i,j,l\in\{1,\dots,d\}$ and $\pmb{x}\in U$,
\Eq{*}{
 \partial_i M_{f,g,m;\mu}(\pmb{x})=\partial_i M_{h,k,m;\mu}(\pmb{x}),&\qquad
 \partial_i\partial_j M_{f,g,m;\mu}(\pmb{x})=\partial_i\partial_j M_{h,k,m;\mu}(\pmb{x}), \\
 \partial_i\partial_j\partial_l M_{f,g,m;\mu}(\pmb{x})&=\partial_i\partial_j\partial_l M_{h,k,m;\mu}(\pmb{x}).
}
In particular, for all $i,j,l\in\{1,\dots,d\}$ and $x\in D$, we have
\Eq{fh}{
 \partial_i M_{f,g,m;\mu}(x^{(d)})=\partial_i M_{h,k,m;\mu}(x^{(d)}),&\qquad
 \partial_i\partial_j M_{f,g,m;\mu}(x^{(d)})=\partial_i\partial_j M_{h,k,m;\mu}(x^{(d)}), \\
 \partial_i\partial_j\partial_l M_{f,g,m;\mu}(x^{(d)})&=\partial_i\partial_j\partial_l M_{h,k,m;\mu}(x^{(d)}).
}
In order to show that the two identities in \eq{1.5} hold on $I$, let $x\in D$ be fixed.
Inequality \eq{MT-1} implies that, for some $i\in\{1,\dots,d\}$,
\Eq{*}{
\mu\Big(\big\{t\in T\mid\partial_i^*m\big(x^{(d)},t\big)\neq0\big\}\Big)>0
}
Therefore, there exists a set $S\subseteq T$ of positive $\mu$-measure such that
$\big(\partial_i^*m\big(x^{(d)},t\big)\big)^2>0$ holds for all $t\in S$. This yields that
\Eq{pos}{
  \big\langle (\partial_i^*m)^2\big\rangle_\mu\big(x^{(d)}\big)>0.
}
Using the second equality in \eq{fh} for $j=i$, and applying formula \eq{PD2}, we get that
\Eq{*}{
  \Phi_{f,g}(x)\big\langle \partial_i^*m)^2\big\rangle_\mu\big(x^{(d)}\big)
  =\Phi_{h,k}(x)\big\langle \partial_i^*m)^2\big\rangle_\mu\big(x^{(d)}\big).
}
In view of \eq{pos}, this equality implies that, for all $x\in D$,
\Eq{Phi}{
\Phi_{f,g}(x)=\Phi_{h,k}(x).
}
By the density of $D$ in $I$ and the continuity of the functions $\Phi_{f,g}$ and $\Phi_{h,k}$, we obtain
that these functions are identical on $I$. Hence the first equality in \eq{1.5} has been verified.

Observe that, until now, we have used only twice continuous differentiability assumptions. The third-order
differentiability will only be used to derive the second equality in \eq{1.5}.

By the assumptions of the theorem, for $x\in D$, there exists $i,j,l\in \{1,\dots,d\}$ such that \eq{MT0}
holds. The third equality in \eq{fh} combined with formula \eq{PD3}, and then the identity
$\Phi_{f,g}=\Phi_{h,k}$ now imply that
\Eq{*}{
  \Psi_{f,g}(x)\big\langle\partial_i^*m\,\partial_j^*m\,\partial_l^*m\big\rangle_\mu\big(x^{(d)}\big)
  =\Psi_{h,k}(x)\big\langle\partial_i^*m\,\partial_j^*m\,\partial_l^*m\big\rangle_\mu\big(x^{(d)}\big).
}
Using condition \eq{MT0}, for all $x\in D$, this simplifies to
\Eq{Psi}{
\Psi_{f,g}(x)=\Psi_{h,k}(x).
}
The density of $D$ in $I$ and the continuity of the functions $\Psi_{f,g}$ and $\Psi_{h,k}$ yields that these
functions are identical on $I$. Therefore, the second equality in \eq{1.5} has also been shown.
\end{proof}

In the next corollary we consider the particular case of \thm{MT} when the measurable family $m$ is given in
the form
\Eq{m}{
  m(\pmb{x},t)=\varphi_1(t)x_1+\cdots+\varphi_d(t)x_d \qquad(\pmb{x}=(x_1,\dots,x_d)\in I^d,\,t\in T).
}
For a $\mu$ integrable function $\varphi:T\to\R$ define $\varphi^*:T\to\R$ by
\Eq{*}{
  \varphi^*(t):=\varphi(t)-\langle \varphi\rangle_\mu.
}

\Cor{MT}{Let $(f,g),(h,k)\in\C_3(I)$, let $\mu$ be a probability measure on the measurable space
$(T,\A)$, let $\varphi_1,\dots,\varphi_d:T\to[0,1]$ $\mu$-measurable functions with
$\varphi_1+\cdots+\varphi_d=1$ and define the measurable family $m:I^d\times T\to\R$ by \eq{m}. Assume that
\Eq{MT-1+}{
  \mu\Big(\big\{t\in T\mid \varphi_1^*(t)=\cdots=\varphi_d^*(t)=0\big\}\Big)<1
}
and there exist $i,j,l\in\{1,\dots,d\}$ such that
\Eq{MT0+}{
\big\langle\varphi_i^*\,\varphi_j^*\,\varphi_l^*\big\rangle_\mu\neq0.
}
Then the following assertions are equivalent:
\begin{enumerate}[(i)]
 \item For all $\pmb{x}\in I^d$, the equality \eq{MT1} holds;
 \item There exists a dense subset $D\subseteq I$ and an open set $U\subseteq I^d$ containing the subdiagonal
$\{x^{(d)}\mid x\in D\}$ such that, for all $\pmb{x}\in U$, the equality \eq{MT1} holds;
 \item The pairs $(f,g)$ and $(h,k)$ are equivalent.
\end{enumerate}}

\begin{proof} The measurable family $m:I^d\times T\to\R$ is given by \eq{m}, hence $m\in\C_3(I^d\times T)$
and, for all $(\pmb{x},t)\in I^d\times T$ and $i,j,l\in\{1,\dots,d\}$, we have
\Eq{*}{
  \partial_i m(\pmb{x},t)=\varphi_i(t) \qquad\mbox{and}\qquad \partial_i^* m(\pmb{x},t)=\varphi_i^*(t).
}
Therefore, conditions \eq{MT-1+} and \eq{MT0+} are equivalent to \eq{MT-1} and \eq{MT0}, respectively.
Thus, the result is a direct consequence of \thm{MT}.
\end{proof}

The next corollary concerns the case when $T=[0,1]$ and $\mu$ is a probability measure on the sigma algebra of
Borel subsets of $[0,1]$. In this setting, define $\hat{\mu}_1$ to be the \textit{first moment} and $\mu_n$
to be the \textit{$n$th centralized moment} of the measure $\mu$ by
\Eq{*}{
  \hat{\mu}_1:=\int_{[0,1]} t\d\mu(t),\qquad \mu_n:=\int_{[0,1]} (t-\hat{\mu}_1)^n\d\mu(t) \qquad(n\in\N).
}

\Cor{MTS}{Let $(f,g),(h,k)\in\C_3(I)$ such that $g$ and $k$ do not vanish on $I$. Let $\mu$ be a probability
measure on the sigma algebra of Borel subsets of $[0,1]$ with $\mu_2\neq0$ and $\mu_3\neq0$. Then the
following assertions are equivalent:
\begin{enumerate}[(i)]
 \item For all $(x,y)\in I^2$, the equality
 \Eq{MT1S}{
 \left(\frac{f}{g}\right)^{-1}\left(
   \frac{\int_{[0,1]} f\big(tx+(1-t)y\big)\d\mu(t)}{\int_{[0,1]} g\big(tx+(1-t)y)\big)\d\mu(t)}\right)
 =\left(\frac{h}{k}\right)^{-1}\left(
   \frac{\int_{[0,1]} h\big(tx+(1-t)y\big)\d\mu(t)}{\int_{[0,1]} k\big(tx+(1-t)y)\big)\d\mu(t)}\right)
 }
 holds;
 \item There exists a dense subset $D\subseteq I$ and an open set $U\subseteq I^2$ containing the subdiagonal
$\{(x,x)\mid x\in D\}$ such that, for all $(x,y)\in U$, the equality \eq{MT1S} holds;
 \item The pairs $(f,g)$ and $(h,k)$ are equivalent.
\end{enumerate}}

\begin{proof}
For the proof of the result, we will apply \cor{MT} in the case when $d=2$ and the measurable family of means
$m:I^2\times [0,1]\to\R$ is given by
\Eq{mS}{
  m((x,y),t):=tx+(1-t)y \qquad (x,y\in I,\,t\in[0,1]),
}
that is, when $\varphi_1(t):=t$ and $\varphi_2(t):=1-t$ for $t\in[0,1]$. In this case, we have that
$\varphi_1^*(t)=-\varphi_2^*(t)=t-\hat{\mu}_1$ and, it is also immediately seen that conditions \eq{MT-1+}
and \eq{MT0+} are equivalent to the inequalities $\mu_2\neq0$ and $\mu_3\neq0$, respectively. Therefore, the
result directly follows from \cor{MT}.
\end{proof}

The next corollary concerns the equality of nonsymmetric weighted two-variable Bajraktarevi\'c means.

\Cor{MTS+}{Let $(f,g),(h,k)\in\C_3(I)$ such that $g$ and $k$ do not vanish on $I$. Let 
$s\in(0,\frac12)\cup(\frac12,1)$. Then the following assertions are equivalent:
\begin{enumerate}[(i)]
 \item For all $(x,y)\in I^2$, the equality
 \Eq{MT1S+}{
 \left(\frac{f}{g}\right)^{-1}\left(\frac{sf(x)+(1-s)f(y)}{sg(x)+(1-s)g(y)}\right)
 =\left(\frac{h}{k}\right)^{-1}\left(\frac{sh(x)+(1-s)h(y)}{sk(x)+(1-s)k(y)}\right)
 }
 holds;
 \item There exists a dense subset $D\subseteq I$ and an open set $U\subseteq I^2$ containing the subdiagonal
$\{(x,x)\mid x\in D\}$ such that, for all $(x,y)\in U$, the equality \eq{MT1S+} holds;
 \item The pairs $(f,g)$ and $(h,k)$ are equivalent.
\end{enumerate}}

\begin{proof}
Let $s\in(0,\frac12)\cup(\frac12,1)$ and apply the previous corollary for the measure 
$\mu:=(1-s)\delta_0+s\delta_1$. Then $\mu$ is a probability measure on the sigma algebra of the Borels 
sets of $[0,1]$ and, for any continuous function $\varphi:[0,1]\to\R$, we have
\Eq{f1}{
  \int_T \varphi(t)\d \mu(t)=\int_T \varphi(t)\d((1-s)\delta_0+s\delta_1)(t)
  =(1-s)\cdot\varphi(0)+s\cdot\varphi(1).
}
Therefore,
\Eq{*}{
  \hat\mu_1=\int_T t\d\mu(t)=(1-s)\cdot0+s\cdot1=s
}
and 
\Eq{*}{
  \mu_2=\int_T (t-\hat\mu_1)^2\d\mu(t)
       &=(1-s)\cdot(-s)^2+s\cdot(1-s)^2=s(1-s)\neq0, \\
  \mu_3=\int_T (t-\hat\mu_1)^3\d\mu(t)
       &=(1-s)\cdot(-s)^3+s\cdot(1-s)^3=s(1-s)(1-2s)\neq0.
}
Thus $\mu$ possesses the properties required in the previous corollary. To complete the proof of the 
corollary, observe that, for any continuous function $\psi:I\to\R$ and $x,y\in I$, we have
\Eq{*}{
  \int_{[0,1]} \psi\big(tx+(1-t)y\big)\d\mu(t)=s\psi(x)+(1-s)\psi(y)
}
if \eq{f1} is applied to the function $\varphi(t):=\psi(tx+(1-t)y)$. In view of the above equality for 
$\psi\in\{f,g,h,k\}$, equation \eq{MT1S+} is equivalent to \eq{MT1S} and hence \cor{MTS} directly implies 
\cor{MTS+}.
\end{proof}

\section{Equality of generalized quasi-arithmetic means}\setcounter{equation}{0} 

In the following results, we are going to characterize the equality of generalized quasi-arithmetic means in
various settings.

\Thm{MT2}{Let $f,g:I\to\R$ be twice continuously differentiable functions such that $f'$ and $g'$ do not
vanish on $I$. Let $m\in\C_2(I^d\times T)$ be a measurable family of means, and let $\mu$ be a probability
measure on the measurable space $(T,\A)$. Assume that, there exists a dense subset $D\subseteq I$ such that,
for all $x\in D$, condition \eq{MT-1} holds. Then the following assertions are equivalent:
\begin{enumerate}[(i)]
 \item For all $\pmb{x}\in I^d$,
 \Eq{MT1+}{
   f^{-1}\left(\int_T f\big(m(\pmb{x},t)\big)\d\mu(t)\right)
   =g^{-1}\left(\int_T g\big(m(\pmb{x},t)\big)\d\mu(t)\right);
 }
 \item There exists an open set $U\subseteq I^d$ containing the subdiagonal $\{x^{(d)}\mid x\in D\}$
such that, for all $\pmb{x}\in U$, the equality \eq{MT1+} holds;
 \item The functions $f''/f'$ and $g''/g'$ are identical on $I$;
 \item There exist real constants $a,b$ such that $g=af+b$.
\end{enumerate}}

\begin{proof}
The implications \textit{(i)}$\Rightarrow$\textit{(ii)} and \textit{(iv)}$\Rightarrow$\textit{(i)} are
trivial. The implication \textit{(iii)}$\Rightarrow$\textit{(iv)} can be seen directly by integrating the
equality $f''/f'=g''/g'$ twice. Therefore, it remains to show that condition \textit{(ii)} implies
\textit{(iii)}.

First observe that the regularity conditions imply that $(f,1),(g,1)\in\C_2(I)$. We can also see that
\eq{MT1+} is equivalent to the equality of the two generalized Bajraktarevi\'c means $M_{f,1,m;\mu}$ and
$M_{g,1,m;\mu}$. Now repeating the same argument that was followed in the proof of \thm{MT}, we can deduce
(under twice differentiability assumptions), that $\Phi_{f,1}$ is equal to $\Phi_{g,1}$ on $I$. This yields
that $f''/f'$ and $g''/g'$ are identical on $I$.
\end{proof}

\Cor{MT2}{Let $f,g:I\to\R$ be twice continuously differentiable functions such that $f'$ and $g'$ do not
vanish on $I$, let $\mu$ be a probability measure on the measurable space $(T,\A)$, let
$\varphi_1,\dots,\varphi_d:T\to[0,1]$ $\mu$-measurable functions with $\varphi_1+\cdots+\varphi_d=1$  such
that condition \eq{MT-1+} holds. Then the following assertions are equivalent:
\begin{enumerate}[(i)]
 \item For all $(x_1,\dots,x_d)\in I^d$,
 \Eq{MT1++}{
   f^{-1}\left(\int_T f\big(\varphi_1(t)x_1+\cdots+\varphi_d(t)x_d)\big)\d\mu(t)\right)
   =g^{-1}\left(\int_T g\big(\varphi_1(t)x_1+\cdots+\varphi_d(t)x_d\big)\d\mu(t)\right);
 }
 \item There exists an open set $U\subseteq I^d$ containing the subdiagonal $\{x^{(d)}\mid x\in D\}$
such that, for all $(x_1,\dots,x_d)\in U$, the equality \eq{MT1++} holds;
 \item There exist real constants $a,b$ such that $g=af+b$.
\end{enumerate}
}

The proof of this corollary is based on \thm{MT2} and can be elaborated in the same way as the proof of
\cor{MT}, the details are left to the reader.

The following consequence of \cor{MT2} has been dealt with in the paper \cite[Theorem 7]{MakPal08}. There 
$f$ and $g$ are assumed only to be continuous, however, the equivalence to condition (ii) is missing.

\Cor{MT2S}{Let $f,g:I\to\R$ be twice continuously differentiable functions such that $f'$ and $g'$ do not
vanish on $I$. Let $\mu$ be a probability measure on the sigma algebra of Borel subsets of $[0,1]$ with
$\mu_2\neq0$. Then the following assertions are equivalent:
\begin{enumerate}[(i)]
 \item For all $(x,y)\in I^2$, the equality
  \Eq{MT1++S}{
   f^{-1}\left(\int_T f\big(tx+(1-t)y\big)\d\mu(t)\right)
   =g^{-1}\left(\int_T g\big(tx+(1-t)y\big)\d\mu(t)\right);
 }
 \item There exists a dense subset $D\subseteq I$ and an open set $U\subseteq I^2$ containing the subdiagonal
$\{(x,x)\mid x\in D\}$ such that, for all $(x,y)\in U$, the equality \eq{MT1++S} holds;
 \item There exist real constants $a,b$ such that $g=af+b$.
\end{enumerate}}

This result follows exactly in the same way from \cor{MT2} as \cor{MTS} follows from \cor{MT}. The next 
statement is related to the equality problem of weighted two-variable quasi-arithmetic means. We note that 
the equivalence of conditions (i) and (iii) can be obtained under the assumption of continuity of the 
generating functions $f$ and $g$. For further and important particular cases of \cor{MT2S}, we refer to the 
examples elaborated in the paper \cite{MakPal08}.

\Cor{MT2S+}{Let $f,g:I\to\R$ be twice continuously differentiable functions such that $f'$ and $g'$ do not
vanish on $I$. Let $s\in(0,1)$. Then the following assertions are equivalent:
\begin{enumerate}[(i)]
 \item For all $(x,y)\in I^2$, the equality
  \Eq{MT1++S+}{
   f^{-1}(sf(x)+(1-s)f(y))=g^{-1}(sg(x)+(1-s)g(y));
 }
 \item There exists a dense subset $D\subseteq I$ and an open set $U\subseteq I^2$ containing the subdiagonal
$\{(x,x)\mid x\in D\}$ such that, for all $(x,y)\in U$, the equality \eq{MT1++S+} holds;
 \item There exist real constants $a,b$ such that $g=af+b$.
\end{enumerate}}

The proof of the above corollary is analogous to that of \cor{MTS+}. It can be deduced from \cor{MT2S} by 
taking the measure $\mu:=(1-s)\delta_0+s\delta_1$ and observing that $\mu_2=s(1-s)\not=0$.

\section{Homogeneity of generalized Bajraktarevi\'c means}\setcounter{equation}{0} 

In this section we characterize the homogeneity of generalized Bajraktarevi\'c and quasi-arithmetic means
under 3 times and 2 times differentiability assumptions, respectively.

Given a nonempty open subinterval $I$ of $\R_+$ and $c>0$, introduce the following notations:
\Eq{*}{
  cI:=\{cx\mid x\in I\}\qquad\mbox{and}\qquad I/I:=\{x/y\mid x,y\in I\}.
}
These sets are also open subintervals of $\R_+$ and the interval $I/I$ is logarithmically symmetric with
respect to $1$, i.e., $u\in I/I$ holds if and only if $1/u\in I/I$. It is also easy to see that the
intersection $I_\lambda:=I\cap\big(\frac{1}{\lambda}I\big)$ is nonempty if and only if $\lambda\in I/I$.

A $d$-variable mean $M:I^d\to\R$ is called \textit{homogeneous} if, for all $\lambda\in I/I$ and for all
$\pmb{x}\in I_\lambda^d$,
\Eq{*}{
  M(\lambda\pmb{x})=\lambda M(\pmb{x}).
}
We will also use the following notation: For a function $f:I\to\R$ and number $\lambda>0$, the function
$f_\lambda:\big(\frac{1}{\lambda}I\big)\to\R$ is defined by
\Eq{*}{
  f_\lambda(x)=f(\lambda x).
}

\Lem{hom}{Assume that $m\in\C_1(I^d\times T)$ is a homogeneous measurable family
of means, and $\mu$ is a probability measure on the measurable space $(T,\A)$. Then, for all 
$i\in\{1,\dots,d\}$ and for all $t\in T$, the mapping
\Eq{*}{
  I\ni x\mapsto\partial_i^*m(x^{(d)},t)
}
is constant on $I$.}

\begin{proof} By the homogeneity of the measurable family $m$, for all $\lambda\in I/I$ and for all
$\pmb{x}\in I_\lambda^d$, we have that 
\Eq{*}{
  m(\lambda\pmb{x},t)=\lambda m(\pmb{x},t).
}
Differentiating this identity with respect to the $i$th variable, we get
\Eq{*}{
  \lambda\partial_i m(\lambda\pmb{x},t)=\lambda\partial_i m(\pmb{x},t), 
}
which simplifies to
\Eq{*}{
  \partial_i m(\lambda\pmb{x},t)=\partial_i m(\pmb{x},t).
}
Let $x,y\in I$ be arbitrary elements. Then, taking $\lambda:=y/x$, and $\pmb{x}:=x^{(d)}$, the above equality 
yields that
\Eq{*}{
  \partial_i m(y^{(d)},t)=\partial_i m(x^{(d)},t)
}
hold for all $x,y\in I$ and $t\in T$. Using this identity, the statement of the lemma follows immediately.
\end{proof}

\Thm{MTH}{Let $(f,g)\in\C_3(I)$, let $m\in\C_3(I^d\times T)$ be a homogeneous measurable family
of means, and let $\mu$ be a probability measure on the measurable space $(T,\A)$. Assume that there exists a
point $x_0\in I$ such that
\Eq{MT-1H}{
  \mu\Big(\big\{t\in T\mid
\partial_1^*m\big(x_0^{(d)},t\big)=\cdots=\partial_d^*m\big(x_0^{(d)},t\big)=0\big\}\Big)<1
}
and there exist $i,j,l\in\{1,\dots,d\}$ such that
\Eq{MT0H}{
\big\langle\partial_i^*m\,\partial_j^*m\,\partial_l^*m\big\rangle_\mu\big(x_0^{(d)}\big)\neq0.
}
Then the following assertions are equivalent:
\begin{enumerate}[(i)]
 \item $M_{f,g,m;\mu}$ is homogeneous;
 \item For all $\lambda\in I/I$ and for all $\pmb{x}\in I_\lambda^d$,
 \Eq{*}{
   M_{f,g,m;\mu}(\pmb{x})=M_{f_\lambda,g_\lambda,m;\mu}(\pmb{x});
 }
 \item For all $\lambda\in I/I$, the pairs $(f,g)$ and $(f_\lambda,g_\lambda)$ are equivalent on the interval
$I_\lambda$;
 \item For all $\lambda\in I/I$ and for all $x\in I_\lambda$,
 \Eq{*}{
   \Phi_{f,g}(x)=\Phi_{f_\lambda,g_\lambda}(x) \qquad\mbox{and}\qquad
   \Psi_{f,g}(x)=\Psi_{f_\lambda,g_\lambda}(x);
 }
 \item There exist two real numbers $\alpha,\beta$ such that $y=f$ and $y=g$ are solutions of the
second-order linear differential equation
 \Eq{Ca}{
   y''(x)=\frac{\alpha}{x}y'(x)+\frac{\beta}{x^2}y(x) \qquad(x\in I);
 }
 \item There exists a pair $(p,q)\in\{(z,w)\in\CC^2\mid z+w,zw\in\R\}$ such that $M_{f,g,m;\mu}$ is
equal to the $d$-variable generalized Gini mean $G_{p,q,m;\mu}$.
\end{enumerate}}

\begin{proof}
Observe first that, based on \lem{hom}, if there exists $x_0\in I$ and $i,j,l\in\{1,\dots,d\}$ such that 
conditions \eq{MT-1H} and \eq{MT0H} hold, then they are also satisfied for all $x_0\in I$.

In order to simplify the computations, in view of \lem{BP}, we may assume that $g$ is
positive, $f/g$ is strictly monotone with a non-vanishing first-order derivative.

Let $M_{f,g,m;\mu}$ be a homogeneous mean. Then, for all $\lambda\in I/I$ and $\pmb{x}\in I_\lambda^d$, we
have
\Eq{*}{
 \frac{1}{\lambda}M_{f_,g_,m;\mu}(\lambda \pmb{x})=M_{f,g,m;\mu}(\pmb{x}).
}
Fix $\lambda\in I/I$ arbitrarily. Using that $m$ is a homogeneous measurable family of means and the
definition of the functions $f_\lambda, g_\lambda$, we can rewrite the last equality as follows
\Eq{lambda}{
 \frac{1}{\lambda}\left(\frac{f}{g}\right)^{-1}\left(
            \frac{\int_T f_\lambda\big(m(\pmb{x},t)\big)\d\mu(t)}
                 {\int_T g_\lambda\big(m(\pmb{x},t)\big)\d\mu(t)}\right)=\left(\frac{f}{g}\right)^{-1}
            \left(\frac{\int_T f\big(m(\pmb{x},t)\big)\d\mu(t)}
                 {\int_T g\big(m(\pmb{x},t)\big)\d\mu(t)}\right).
}
Also, one can easily see that, for all $u\in(f/g)(I)$, we get
\Eq{*}{
  \frac{1}{\lambda} \left(\frac{f}{g}\right)^{-1}(u)=\left(\frac{f_\lambda}{g_\lambda}\right)^{-1}(u),
}
hence the equality in \eq{lambda} reduces to
\Eq{*}{
 M_{f_{\lambda},g_{\lambda},m;\mu}(\pmb{x})=M_{f,g,m;\mu}(\pmb{x}).
}
Therefore, the homogeneity of the mean $M_{f,g,m;\mu}$ implies \textit{(ii)}, that is the equality of the two
$d$-variable generalized Bajraktarevi\'c means $M_{f,g,m;\mu}$ and $M_{f_{\lambda},g_{\lambda},m;\mu}$ on
$I_\lambda^d$. In fact, from this argument, also the equivalence of these statements can be seen.

Applying \thm{MT} for $h:=f_\lambda$, $k:=g_\lambda$, we get that \textit{(ii)},  \textit{(iii)}, and
\textit{(iv)} are equivalent to each other and therefore to \textit{(i)}, too.

In order to understand the content of condition \textit{(iv)}, observe first that, for all $x\in\frac1\lambda
I$,
\Eq{*}{
  \Phi_{f_\lambda,g_\lambda}(x)
  =\frac{\partial_1^2\D_{f_\lambda,g_\lambda}(x,x)}{\partial_1\D_{f_\lambda,g_\lambda}(x,x)}
  =\lambda\frac{\partial_1^2\D_{f,g}(\lambda x,\lambda x)}{\partial_1\D_{f,g}(\lambda x,\lambda x)}
  =\lambda\Phi_{f,g}(\lambda x)
}
and similarly,
\Eq{*}{
  \Psi_{f_\lambda,g_\lambda}(x)
  =\frac{\partial_1^2\partial_2\D_{f_\lambda,g_\lambda}(x,x)}{\partial_1\D_{f_\lambda,g_\lambda}(x,x)}
  =\lambda^2\frac{\partial_1^2\partial_2\D_{f,g}(\lambda x,\lambda x)}
     {\partial_1\D_{f,g}(\lambda x,\lambda x)}
  =\lambda^2\Psi_{f,g}(\lambda x).
}
Therefore, condition \textit{(iv)} holds if and only if, for all $\lambda\in I/I$ and $x\in I_\lambda$,
\Eq{uv}{
  \lambda\Phi_{f,g}(\lambda x)=\Phi_{f,g}(x)
  \qquad\mbox{and}\qquad
  \lambda^2\Psi_{f,g}(\lambda x)=\Psi_{f,g}(x).
}
We now show that the maps $x\mapsto x\Phi_{f,g}(x)$ and $x\mapsto x^2\Psi_{f,g}(x)$ are constants over $I$.
To see this, let $u,v\in I$ be arbitrary. Then $\lambda:=v/u\in I/I$ and $x:=u\in I_\lambda$.
Therefore, by first equality in \eq{uv},
\Eq{*}{
  u\Phi_{f,g}(u)=u\Phi_{f,g}(x)=u\lambda\Phi_{f,g}(\lambda x)=v\Phi_{f,g}(v),
}
which proves that $x\mapsto x\Phi_{f,g}(x)$ is a constant map on $I$. Denoting the value of this map by
$\alpha$, we get that
\Eq{one}{
  \Phi_{f,g}(x)=\frac{\alpha}{x} \qquad(x\in I).
}
A completely similar argument applied for the second identity in \eq{uv} shows that there exists a constant
$\beta\in\R$ such that
\Eq{two}{
  \Psi_{f,g}(x)=\frac{\beta}{x^2} \qquad(x\in I).
}
One can also see that these identities are also sufficient for \eq{uv} to hold.

As we have seen in the proof of \thm{1.5}, the functions $f$ and $g$ are solutions of the second order
homogeneous linear differential equation \eq{y''}, hence, by the formulae \eq{one} and \eq{two}, we can see
that condition \textit{(v)} is fulfilled. To prove the implication \textit{(v)}$\Rightarrow$\textit{(vi)},
assume that $y=f$ and $y=g$ are solutions of the homogeneous second-order Cauchy--Euler equation \eq{Ca} for
some $\alpha,\beta\in\R$. For a function $y:I\to\R$ define $Y:=y\circ\exp$. Then
$y=Y\circ\log$, and using the chain rule, we deduce that
\Eq{*}{
y'(x)=\frac{1}{x}Y'(\log x)\qquad\mbox{and}\qquad y''(x)=\frac{1}{x^2}\big(Y''(\log x)-Y'(\log x)\big).
}
Substituting these expressions into equation \eq{Ca}, we get that $y:I\to\R$ is a solution of \eq{Ca} if and
only if $Y:\log(I)\to\R$ is the solution of the following homogeneous second-order differential equation with
constant coefficients
\Eq{Ca+}{
Y''(u)=(\alpha +1)Y'(u)+\beta Y(u) \qquad(u\in \log(I)).
}
Denote by $p$ and $q$ the roots of the characteristic polynomial
\Eq{ro}{
r^2-(\alpha+1)r-\beta=0.
}
Then, we have that $(p,q)\in\{(z,w)\in\CC^2\mid z+w,zw\in\R\}$. Therefore, either $p,q\in\R$ or
$p=\bar{q}\in\CC\setminus\R$. According to these possibilities, a fundamental system
$Y_1,Y_2:\log(I)\to\R$ for the solution of \eq{Ca+} can be obtained in the following form
\Eq{*}{
  Y_{1}(u)&=\exp(pu),&\quad  Y_{2}(u)&=\exp(qu)& \quad&\mbox{if}\quad p,q\in\R,\,p\neq q,\\
  Y_{1}(u)&=\exp(pu),&\quad  Y_{2}(u)&=u\exp(pu)& \quad&\mbox{if}\quad p=q\in\R,\\
  Y_{1}(u)&=\exp(au)\cos(bu),&\quad  Y_{2}(u)&=\exp(au)\sin(bu)& \quad&\mbox{if}\quad
      p=\bar{q}=a+bi,\,\,a,b\in\R,\,b\neq0.
}
Substituting $u=\log x$, we get that the functions $y_1,y_2:I\to\R$ defined by
\Eq{Sol+}{
  y_1(x)&=x^p,&\quad y_2(x)&=x^q& \qquad&\mbox{if}\quad p,q\in\R,\,p\neq q,\\
  y_1(x)&=x^p\log(x),&\quad y_2(x)&=x^p& \qquad&\mbox{if}\quad p=q\in\R,\\
  y_1(x)&=x^a\cos(b\log(x)),&\quad y_2(x)&=x^a\sin(b\log(x))& \qquad&\mbox{if}\quad
      p=\bar{q}=a+bi,\,\,a,b\in\R,\,b\neq0.
}
form a fundamental system of differential equation \eq{Ca}. Therefore, $f$ and $g$ are linear combinations of
$y_1$ and $y_2$. Thus the pairs $(f,g)$ and $(y_1,y_2)$ determine the same Bajraktarevi\'c mean (provided
that they belong to $\C_3(I)$. It is obvious that these functions are infinitely many times differentiable.
One can also easily see that their Wronski determinant does not vanish anywhere. As we have seen it in the 
introduction, these functions always form a Chebyshev system on $\R_+$ if $p,q\in\R$. On the other hand, if 
$p=\bar{q}\not\in\R$, then $(y_1,y_2)$ is a Chebyshev system over $I$ if and only if $I/I$ is contained in 
the open interval $\big(\exp\big(-\tfrac{\pi}{|b|}\big),\exp\big(\tfrac{\pi}{|b|}\big)\big)$.

This implies that, in each of the above cases, $M_{f,g,m;\mu}$ is equal to the $d$-variable generalized Gini
mean $G_{p,q,m;\mu}$ (which we defined in the introduction). Hence the implication 
\textit{(v)}$\Rightarrow$\textit{(vi)} is satisfied.

To complete the proof of the theorem, we can easily prove the implication
\textit{(vi)}$\Rightarrow$\textit{(iv)}. Assume that, for some $(p,q)\in\{(z,w)\in\CC^2\mid z+w,zw\in\R\}$,
we have
\Eq{*}{
M_{f,g,m;\mu}(\pmb{x})=G_{p,q,m;\mu}(\pmb{x}) \qquad(\pmb{x}\in I^d),
}
then \thm{MT} implies that the pairs $(f,g)$ and $(y_1,y_2)$ are equivalent. Therefore, \thm{1.5} implies
that
\Eq{*}{
\Phi_{f,g}(x)=\Phi_{y_1,y_2}(x)=\frac{p+q-1}{x}\qquad\mbox{and}\qquad 
\Psi_{f,g}(x)=\Psi_{y_1,y_2}(x)=\frac{pq}{x^2}.
}
Then, for all $\lambda\in I/I$ and $x\in I_\lambda$, the two equations of condition \textit{(iv)} are
trivially satisfied. This completes the proof of the theorem.
\end{proof}

In the following result, we derive an important particular case of \thm{MTH} when the measurable family of 
means consists of weighted $d$-variable arithmetic means.

\Cor{MTH1}{Let $(f,g)\in\C_3(I)$, let $\mu$ be a probability measure on the measurable space
$(T,\A)$, let $\varphi_1,\dots,\varphi_d:T\to[0,1]$ $\mu$-measurable functions with
$\varphi_1+\cdots+\varphi_d=1$ and define the measurable family $m:I^d\times T\to\R$ by \eq{m}.
Assume that conditions \eq{MT-1+} and \eq{MT0+} are satisfied. Then $M_{f,g,m;\mu}$ is homogeneous if and only 
if there exists a pair $(p,q)\in\{(z,w)\in\CC^2\mid z+w,zw\in\R\}$ such that $M_{f,g,m;\mu}$ is equal to the 
$d$-variable generalized Gini mean $G_{p,q,m;\mu}$.}

\begin{proof}
Applying the same argument as in the proof of \cor{MT}, the verification of this corollary directly 
follows from the previous theorem because the measurable family $m$ is trivially homogeneous and infinitely 
many times differentiable, furthermore, \eq{MT-1+} and \eq{MT0+} imply conditions \eq{MT-1H} and \eq{MT0H}.
\end{proof}

The following two results are direct consequences of \cor{MTH1}. Their proofs go on the same line as the 
proofs of \cor{MTS} and \cor{MTS+}, respectively but using \cor{MTH1} instead of \cor{MT}.

\Cor{MTSH+}{Let $(f,g)\in\C_3(I)$ such that $g$ does not vanish on $I$. Let $\mu$ be a probability
measure on the sigma algebra of Borel subsets of $[0,1]$ with $\mu_2\neq0$ and $\mu_3\neq0$. Then the
the functional equation 
 \Eq{*}{
 \left(\frac{f}{g}\right)^{-1}\left(
   \frac{\int_{[0,1]} f\big(t\lambda x+(1-t)\lambda y\big)\d\mu(t)}
   {\int_{[0,1]} g\big(t\lambda x+(1-t)\lambda y\big)\d\mu(t)}\right)
 =\lambda\left(\frac{f}{g}\right)^{-1}\left(
   \frac{\int_{[0,1]} f\big(tx+(1-t)y\big)\d\mu(t)}{\int_{[0,1]}g\big(tx+(1-t)y\big)\d\mu(t)}\right)
 }
holds for all $\lambda\in I/I$ and for all $(x,y)\in I^2$ if and only if there exists a pair 
$(p,q)\in\{(z,w)\in\CC^2\mid z+w,zw\in\R\}$ such that, for all $x,y\in I$, the expression
\Eq{*}{
 \left(\frac{f}{g}\right)^{-1}\left(
   \frac{\int_{[0,1]} f\big(tx+(1-t)y\big)\d\mu(t)}{\int_{[0,1]} g\big(tx+(1-t)y\big)\d\mu(t)}\right)
}
is of the form
{\small
\Eq{*}{
      \begin{cases}
          \left(\dfrac{\int_{[0,1]} \big(tx+(1-t)y)\big)^p\d\mu(t)}
                 {\int_{[0,1]} \big(tx+(1-t)y)\big)^q\d\mu(t)}\right)^{\frac1{p-q}} 
                 &\mbox{if }p,q\in\R,\,p\neq q,\\[5mm]
          \exp\left(\dfrac{\int_{[0,1]} \big(tx+(1-t)y)\big)^p\log\big(tx+(1-t)y)\big)\d\mu(t)}
                 {\int_{[0,1]} \big(tx+(1-t)y\big)^p\d\mu(t)}\right) &\mbox{if }p=q\in\R,\\[5mm]
          \exp\left(\dfrac{1}{b}\arctan\left(\dfrac{\int_{[0,1]} 
                 \big(tx+(1-t)y)\big)^a\sin\big(b\log\big(tx+(1-t)y)\big)\big)\d\mu(t)}   
                 {\int_{[0,1]} \big(tx+(1-t)y)\big)^a\cos\big(b\log\big(tx+(1-t)y)\big)\big) 
                    \d\mu(t)}\right)\right) &\mbox{if }p=\bar{q}=a+bi,\,b\neq0,
        \end{cases}
}}
provided that, in the last case, the inclusion
$I\subseteq\big(\exp\big(-\frac{\pi}{2|b|}\big),\exp\big(\frac{\pi}{2|b|}\big)\big)$ holds.
}

The result presented in the above corollary was established by Losonczi in \cite[Theorem 2.1]{Los13} under 
the assumption of six times continuous differentiability of $f,g$ and the moment conditions $\mu_2\neq0$ and 
$5\mu_2\mu_4^2+\mu_4\mu_6-6\mu_2^2\mu_6\neq0$. Having a careful look at the proof of this result in
\cite{Los13}, under the condition $\mu_3\neq0$, the above conclusion was reached using only three 
times differentiability.

\Cor{MTSH++}{Let $(f,g)\in\C_3(I)$ such that $g$ does not vanish on $I$. Let $s\in(0,\frac12)\cup(\frac12,1)$. 
Then the functional equation
 \Eq{*}{
 \left(\frac{f}{g}\right)^{-1}\left(
   \frac{sf(\lambda x)+(1-s)f(\lambda y)}{sg(\lambda x)+(1-s)g(\lambda y)}\right)
   =\lambda\left(\frac{f}{g}\right)^{-1}\left(\frac{sf(x)+(1-s)f(y)}{sg(x)+(1-s)g(y)}\right)
 }
 holds for all $\lambda\in I/I$ and for all $(x,y)\in I^2$ if and only if there exists a pair 
$(p,q)\in\{(z,w)\in\CC^2\mid z+w,zw\in\R\}$ such that, for all $x,y\in I$, the expression
\Eq{*}{
 \left(\frac{f}{g}\right)^{-1}\left(\frac{sf(x)+(1-s)f(y)}{sg(x)+(1-s)g(y)}\right)
}
is of the form
\Eq{*}{
      \begin{cases}
          \left(\dfrac{sx^p+(1-s)y^p}{sx^q+(1-s)y^q}\right)^{\frac1{p-q}} 
                 &\mbox{if }p,q\in\R,\,p\neq q,\\[5mm]
          \exp\left(\dfrac{sx^p\log(x)+(1-s)y^p\log(y)}{sx^p+(1-s)y^p}\right) &\mbox{if }p=q\in\R,\\[5mm]
          \exp\left(\dfrac{1}{b}\arctan\left(\dfrac{sx^a\sin(\log(x^b))+(1-s)y^a\sin(\log(y^b))}   
                 {sx^a\cos(\log(x^b))+(1-s)y^a\cos(\log(y^b))}\right)\right) 
                 &\mbox{if }p=\bar{q}=a+bi,\,b\neq0,
        \end{cases}
}
provided that, in the last case, the inclusion
$I\subseteq\big(\exp\big(-\frac{\pi}{2|b|}\big),\exp\big(\frac{\pi}{2|b|}\big)\big)$ holds.
}

\section{Homogeneity of generalized quasi-arithmetic means}\setcounter{equation}{0} 

In this section we consider particular cases of the results of Section 5, when all the means are generalized 
quasi-arithmetic. Let $I$ again be an open subinterval of $\R_+$. The proofs of the results can be obtained 
by combining the arguments of Sections 4 and 5 therefore they are completely omitted.

\Thm{MTH2}{Let $f:I\to\R$ be a twice continuously differentiable function such that $f'$ does not
vanish on $I$. Let $m\in\C_2(I^d\times T)$ be a homogeneous measurable family
of means, and let $\mu$ be a probability measure on the measurable space $(T,\A)$. 
Assume that condition \eq{MT-1H} holds. Then the following assertions are equivalent:
\begin{enumerate}[(i)]
 \item $M_{f,1,m;\mu}$ is homogeneous;
 \item For all $\lambda\in I/I$ and for all $\pmb{x}\in I_\lambda^d$,
 \Eq{*}{
   M_{f,1,m;\mu}(\pmb{x})=M_{f_\lambda,1,m;\mu}(\pmb{x});
 }
 \item For all $\lambda\in I/I$, there exist real constants $a_\lambda,b_\lambda$ such that 
$f_\lambda(x)=a_\lambda f(x)+b_\lambda$ holds for all $x\in I_\lambda$;
 \item For all $\lambda\in I/I$, the functions $f''/f'$ and $f''_\lambda/f'_\lambda$ 
 are identical on $I_\lambda$,
 \item There exists a real number $\alpha$ such that $y=f$ is a solution of the
second-order linear differential equation
 \Eq{*}{
   y''(x)=\frac{\alpha}{x}y'(x) \qquad(x\in I);
 }
 \item There exists a real number $p$ such that $M_{f,1,m;\mu}$ is equal to the $d$-variable 
 generalized H\"older mean $H_{p,m;\mu}$.
\end{enumerate}}

The above theorem reduces to the following result if the measurable family consists of weighted arithmetic 
means.

\Cor{MTH2}{Let $f:I\to\R$ be a twice continuously differentiable function such that $f'$ does not
vanish on $I$. Let $\mu$ be a probability measure on the measurable space $(T,\A)$, let 
$\varphi_1,\dots,\varphi_d:T\to[0,1]$ be $\mu$-measurable functions with $\varphi_1+\cdots+\varphi_d=1$ and 
define the measurable family $m:I^d\times T\to\R$ by \eq{m}. Assume that conditions \eq{MT-1+} is satisfied. 
Then the $d$-variable generalized quasi-arithmetic mean $M_{f,1,m;\mu}$ is homogeneous if and only if there 
exists a real number $p$ such that $M_{f,1,m;\mu}$ is equal to the $d$-variable generalized H\"older mean  
$H_{p,m;\mu}$.}

The following particular case of \cor{MTH2} was obtained by Burai and Jarczyk \cite{BurJar13} in 2013. 

\Cor{MTSH+q}{Let $f:I\to\R$ be a twice continuously differentiable function such that $f'$ does not
vanish on $I$. Let $\mu$ be a probability measure on the sigma algebra of Borel subsets of $[0,1]$ with 
$\mu_2\neq0$ and $\mu_3\neq0$. Then the functional equation 
 \Eq{*}{
 f^{-1}\bigg(
   \int\limits_{[0,1]} f\big(t\lambda x+(1-t)\lambda y\big)\d\mu(t)\bigg)
 =\lambda f^{-1}\bigg(\int\limits_{[0,1]} f\big(tx+(1-t)y\big)\d\mu(t)\bigg)
 }
holds for all $\lambda\in I/I$ and for all $(x,y)\in I^2$ if and only if there exists $p\in\R$ such that, for 
all $x,y\in I$,
\Eq{*}{
 f^{-1}\bigg(\int\limits_{[0,1]}\!\! f\big(tx+(1-t)y\big)\d\mu(t)\bigg)
    =   \begin{cases}
          \bigg(\int\limits_{[0,1]}\!\! \big(tx+(1-t)y)\big)^p\d\mu(t)\bigg)^{\frac1{p}} 
                 &\mbox{if }p\neq 0,\\[5mm]
          \exp\bigg(\int\limits_{[0,1]}\!\! \big(tx+(1-t)y)\big)^p\log\big(tx+(1-t)y)\big)\d\mu(t)\bigg) 
                 &\mbox{if }p=0.
        \end{cases}
}}

Upon taking the particular measure $\mu:=(1-s)\delta_0+s\delta_1$ in the above corollary, we can deduce a 
classical result for the homogeneity of two-variable weighted quasi-arithmetic means (cf.\ 
\cite{HarLitPol34}).

\Cor{MTSH++q}{Let $f:I\to\R$ be a twice continuously differentiable function such that $f'$ does not
vanish on $I$. Let $s\in(0,\frac12)\cup(\frac12,1)$. Then the functional equation
 \Eq{*}{
 f^{-1}(sf(\lambda x)+(1-s)f(\lambda y))
   =\lambda f^{-1}(sf(x)+(1-s)f(y)
 }
holds for all $\lambda\in I/I$ and for all $(x,y)\in I^2$ if and only if there exists a pair 
$p\in\R$ such that, for all $x,y\in I$,
\Eq{*}{
 f^{-1}(sf(x)+(1-s)f(y))
   =   \begin{cases}
          (sx^p+(1-s)y^p)^{\frac1{p}} &\mbox{if }p\neq 0,\\[2mm]
          x^sy^{1-s} &\mbox{if }p=0.
        \end{cases}
}}

\def\MR{}

\begin{thebibliography}{10}

\bibitem{AczDar63c}
{\sc Aczél, J., and Daróczy, Z.}
\newblock {Über verallgemeinerte quasilineare {M}ittelwerte, die mit
  {G}ewichtsfunktionen gebildet sind}.
\newblock {\em Publ. Math. Debrecen 10\/} (1963), 171–190.

\bibitem{Baj58}
{\sc Bajraktarević, M.}
\newblock {Sur une équation fonctionnelle aux valeurs moyennes}.
\newblock {\em Glasnik Mat.-Fiz. Astronom. Društvo Mat. Fiz. Hrvatske Ser. II
  13\/} (1958), 243–248.

\bibitem{Baj69}
{\sc Bajraktarević, M.}
\newblock {Über die {V}ergleichbarkeit der mit {G}ewichtsfunktionen gebildeten
  {M}ittelwerte}.
\newblock {\em Studia Sci. Math. Hungar. 4\/} (1969), 3–8.

\bibitem{BerMor98}
{\sc Berrone, L.~R., and Moro, J.}
\newblock {Lagrangian means}.
\newblock {\em Aequationes Math. 55}, 3 (1998), 217–226.

\bibitem{BerMor00}
{\sc Berrone, L.~R., and Moro, J.}
\newblock {On means generated through the {C}auchy mean value theorem}.
\newblock {\em Aequationes Math. 60}, 1-2 (2000), 1–14.

\bibitem{BesPal03}
{\sc Bessenyei, M., and Páles, Z.}
\newblock {Hadamard-type inequalities for generalized convex functions}.
\newblock {\em Math. Inequal. Appl. 6}, 3 (2003), 379–392.

\bibitem{BurJar13}
{\sc Burai, P., and Jarczyk, J.}
\newblock {Conditional homogeneity and translativity of {M}akó–{P}áles
  means}.
\newblock {\em Annales Univ. Sci. Budapest. Sect. Comp. 40\/} (2013),
  159–172.

\bibitem{Gin38}
{\sc Gini, C.}
\newblock {{D}i una formula compressiva delle medie}.
\newblock {\em Metron 13\/} (1938), 3–22.

\bibitem{HarLitPol34}
{\sc Hardy, G.~H., Littlewood, J.~E., and Pólya, G.}
\newblock {\em {Inequalities}}.
\newblock Cambridge University Press, Cambridge, 1934.
\newblock (first edition), 1952 (second edition).

\bibitem{Los99}
{\sc Losonczi, L.}
\newblock {Equality of two variable weighted means: reduction to differential
  equations}.
\newblock {\em Aequationes Math. 58}, 3 (1999), 223–241.

\bibitem{Los02d}
{\sc Losonczi, L.}
\newblock {Comparison and subhomogeneity of integral means}.
\newblock {\em Math. Inequal. Appl. 5}, 4 (2002), 609–618.

\bibitem{Los02a}
{\sc Losonczi, L.}
\newblock {Homogeneous {C}auchy mean values}.
\newblock In {\em {Functional Equations — Results and Advances}}, Z.~Daróczy
  and Z.~Páles, Eds., vol.~3 of {\em {Advances in Mathematics}}. Kluwer Acad.
  Publ., Dordrecht, 2002, p.~209–218.

\bibitem{Los02c}
{\sc Losonczi, L.}
\newblock {On the comparison of {C}auchy mean values}.
\newblock {\em J. Inequal. Appl. 7}, 1 (2002), 11–24.

\bibitem{Los03a}
{\sc Losonczi, L.}
\newblock {Equality of two variable {C}auchy mean values}.
\newblock {\em Aequationes Math. 65}, 1-2 (2003), 61–81.

\bibitem{Los06b}
{\sc Losonczi, L.}
\newblock {Equality of two variable means revisited}.
\newblock {\em Aequationes Math. 71}, 3 (2006), 228–245.

\bibitem{Los07a}
{\sc Losonczi, L.}
\newblock {Homogeneous non-symmetric means of two variables}.
\newblock {\em Demonstratio Math. 40}, 1 (2007), 169–180.

\bibitem{Los07b}
{\sc Losonczi, L.}
\newblock {Homogeneous symmetric means of two variables}.
\newblock {\em Aequationes Math. 74}, 3 (2007), 262–281.

\bibitem{Los13}
{\sc Losonczi, L.}
\newblock {On homogenous {P}áles means}.
\newblock {\em Ann. Univ. Sci. Budapest. Sect. Comput. 41\/} (2013), 103–117.

\bibitem{LosPal08}
{\sc Losonczi, L., and Páles, Z.}
\newblock {Comparison of means generated by two functions and a measure}.
\newblock {\em J. Math. Anal. Appl. 345}, 1 (2008), 135–146.

\bibitem{LosPal11a}
{\sc Losonczi, L., and Páles, Z.}
\newblock {Equality of two-variable functional means generated by different
  measures}.
\newblock {\em Aequationes Math. 81}, 1-2 (2011), 31–53.

\bibitem{MakPal08}
{\sc Makó, Z., and Páles, Z.}
\newblock {On the equality of generalized quasiarithmetic means}.
\newblock {\em Publ. Math. Debrecen 72\/} (2008), 407–440.

\bibitem{Pal89c}
{\sc Páles, Z.}
\newblock {On comparison of homogeneous means}.
\newblock {\em Ann. Univ. Sci. Budapest. Eötvös Sect. Math. 32\/} (1989),
  261–266 (1990).

\bibitem{Pal11}
{\sc Páles, Z.}
\newblock {On the equality of quasi-arithmetic and {L}agrangian means}.
\newblock {\em J. Math. Anal. Appl. 382}, 1 (2011), 86–96.

\bibitem{PalZak17}
{\sc Páles, Z., and Zakaria, A.}
\newblock {On the local and global comparison of generalized {B}ajraktarević
  means}.
\newblock {\em J. Math. Anal. Appl. 455}, 1 (2017), 792–815.

\end{thebibliography}

\end{document}